\documentclass[12pt]{article}
\usepackage{amsmath,amssymb,amsfonts,amsthm}
\usepackage{eucal}

\def\be{\begin{eqnarray}}
\def\ee{\end{eqnarray}}

\DeclareMathOperator{\Stab}{Stab}

\def\fC{{\frak C}}

\begin{document}
\title{On the Instanton $R$-matrix}

\vspace{3mm}

\author{{ {{Andrey Smirnov}} }\footnote{ {\small {
Department of Mathematics, Columbia University, New York USA} and { ITEP, Moscow, Russia}};
asmirnov@math.columbia.edu, asmirnov@itep.ru}}

%{Andrey Smirnov }
\maketitle

\begin{abstract}
A torus action on a symplectic variety allows one to construct solutions to the quantum Yang-Baxter equations ($R$-matrices).
For a torus action on cotangent bundles over flag varieties the resulting $R$-matrices are the standard rational solutions
of the Yang-Baxter equation, which are well known in the theory of quantum integrable systems. The torus action on the instanton moduli space leads to more complicated $R$-matrices, depending additionally on two equivariant parameters $t_1$ and $t_2$.
In this paper we derive an explicit expression for the $R$-matrix associated with the instanton moduli space. We study its matrix elements and its Taylor expansion in the powers of the spectral parameter. Certain matrix elements of this $R$-matrix give a generating function for the characteristic classes of tautological bundles over the Hilbert schemes in terms of the bosonic cut-and-join operators. In particular we rederive from the $R$-matrix the well known Lehn's formula for the first Chern class.

We explicitly compute  the first several coefficients for the power series expansion of the $R$-matrix in the spectral parameter. These coefficients are represented by simple contour integrals of some symmetrized bosonic fields.
\end{abstract}

\newpage

\tableofcontents

\section{Introduction}
 Many achievements in modern mathematical physics are deeply interrelated and in the end should appear as faces of one mathematical structure.
Among the subjects that, presumably, have common origin are the recent geometrical research of instanton moduli spaces \cite{OM}-\cite{MMN2}, the topological string theories \cite{StrThFirst}-\cite{StrThLast}, supersymmetric quantum field theories \cite{SUSY1}-\cite{SUSY5}, work on quantum invariants of three dimensional manifolds and knots, three-dimensional gravity \cite{Knots1}-\cite{Knots7}, exactly solvable models of classical and statistical mechanics \cite{QISMfirst}-\cite{QISMlast} and other topics.

These theories are distingushed by their \textit{integrability}. This means that many properties of their models can be understood completely using only a big \textit{internal symmetry}  of these models. A good example of this property is the explicit computation of the celebrated Nekrasov partition functions \cite{SUSY1}-\cite{SUSY4}, or recent computations of knot invariants in the refined Chern-Simons theory \cite{Knots2}-\cite{Knots5}.

The symmetries of these models are not necessarily classical objects. In the simplest case of classical integrable systems they are well understood universal enveloping algebras of the simple Lie algebras. In the case of the quantum theories they might form much more complicated structures - Yangians, quantum groups, Hecke algebras, etc. The representation theory for these algebras is becoming one of the main tools in this field \cite{OM}.

The main approach to investigating these algebras was developed in the study of the exactly solvable statistical systems and is known as the \textit{quantum inverse scattering method} (QISM) \cite{QISMfirst}-\cite{MJD}. The main ingredient of QISM is the quantum $R$-matrix. It miraculously contains all the essential information about the symmetry algebra. Once we know the corresponding $R$-matrix, the QISM usually gives an explicit recipe for computing all the essential properties of the system. For example we get correlation functions in TQFT, full set of commuting hamiltonians in the quantum integrable system, invariants of 3d manifolds and knots, the explicit expressions for characteristic classes in $K$-theory of $T$-spaces and more.

Even though the common integrable nature of the fields mentioned above is known, nobody has attempted, to the best of our knowledge, to find the $R$-matrix related to these theories. This article is a first step in this direction.

Following recent work of A. Okounkov and D. Maulik (\cite{OM} section 4), we define the $R$-matrix as the operator acting in $H^{\bullet}_{T}(X)$, where $X$ is a symplectic variety endowed with the hamiltonian action of a torus $T$. By definition, this operator gives a certain solution of the quantum Yang-Baxter equation with a spectral parameter $u\in H^{\bullet}_{T}(\textrm{pt})$. We compute the $R$-matrix explicitly  for the cotangent bundles to Flag varieties, more generally - $A_{n}$ quiver varieties and for $X={\cal{M}}(r,n)$ - the moduli space of rank $r$ framed sheaves on $\mathbb{P}^2$ or simply the instanton moduli space.  The first two examples show that this method reproduces the well known operators in the QISM: for the flag manifolds the corresponding $R$-matrices coincide with the standard rational $R$-matrix for $\frak{gl}(N)$ in the fundamental representations.

 We are mainly interested in the operator ${\cal{R}}(u)$ acting on the cohomology of the instanton moduli space. In section \ref{comp} we derive an explicit expression for this operator in terms of the fermionic operators. In addition, we check that it reproduces known formulae for the equivariant characteristic classes of tautological bundle over Hilbert schemes. For example, using the explicit form of $R$-matrix we rederive Lehn's formula for the first Chern class of the tautological bundle \cite{Lehn}.

We continue with a short outline of the geometry of Hilbert schemes, moduli spaces and $R$-matrices.
The Hilbert scheme $\textrm{Hilb}_n$ of $n$ points in the plane ${\mathbb{C}}^2$ has been and remains the focus of numerous studies in the field of mathematical physics \cite{Ok1}-\cite{Lehn2}.
The cohomology of this space is identified with the space of symmetric polynomials in an infinite number of variables:
\be
H^{\bullet}\Big(\coprod\limits_{n=0}^{\infty} \textrm{Hilb}_n \Big)\simeq  {\mathbb{C}}[p_1,p_2,...]={\cal{F}}
\ee
The last space is usually referred to as the boson Fock space and will be denoted by ${\cal{F}}$. As a ring ${\cal{F}}$ is generated by  the Chern classes $c_{k}({\cal{V}})$ of the tautological bundle ${\cal{V}} $ over $ \coprod\limits_{n=0}^{\infty} \textrm{Hilb}_n$. The multiplication by the Chern classes gives certain operators in ${\cal{F}}$. For example in \cite{Lehn} the following formula for the first Chern class is derived:
\be
\label{fcc}
c_{1}({\cal{V}})=-\dfrac{1}{2}\sum\limits_{m,n=0}^{\infty}\, n m p_{n+m} \dfrac{\partial^2}{ \partial p_{n} \partial p_{m}} = -\dfrac{1}{2}\sum\limits_{m,n=1}^{\infty} \alpha_{-m-n} \alpha_{n} \alpha_{m}
\ee
where the \textit{boson} operators $\alpha_{n}$ act on the Fock space as in \cite{Nak1}:
\be
\label{bosons}
\alpha_{-n}=p_n, \ \ \ \alpha_{n}= n \dfrac{\partial}{\partial p_n}, \ \ \textrm{for} \ \ \ n>0
\ee
The cohomological degree of the polynomials is $\deg p_{n}=2(n-1)$. Thus, (\ref{fcc}) is a locally nilpotent operator on ${\cal{F}}$ increasing the degree by~2.
The equivariant cohomology of this space is a tensor product:
\be
H^{\bullet}_{B}\Big(\coprod\limits_{n=0}^{\infty} \textrm{Hilb}_n \Big) \simeq  H^{\bullet}\Big(\coprod\limits_{n=0}^{\infty} \textrm{Hilb}_n \Big)  \otimes {\mathbb{C}}[t_1,t_2]
\ee
Here $B \simeq ({\mathbb{C}}^{\ast})^2$ is an algebraic torus acting on a plane by scaling the coordinates:
\be
\label{pa}
(z_1,z_2) \cdot(x,y) =(z_1 x, z_2 y)
\ee
what induces an action of $B$ on $\textrm{Hilb}_n$ and $t_1$, $t_2$ are the equivariant parameters corresponding to the characters of $B$.

In the equivariant case, the operators of multiplication on the characteristic classes are not locally nilpotent anymore. The operator for $c_{1}({\cal{V}})$ in the equivariant case has the form \cite{Ok1}:
\begin{small}
\be
\label{lf}
\begin{array}{r}
c_{1}({\cal{V}})=\dfrac{1}{2}\sum\limits_{m,n=1}^{\infty}\Big( t_1 t_2 \alpha_{-m} \alpha_{-n} \alpha_{n+m} - \alpha_{-m-n} \alpha_{n} \alpha_{m}\Big)+
\dfrac{\hbar}{2}\, \sum\limits_{n=1}^{\infty}\,( n-1) \alpha_{-n}\alpha_{n}
\end{array}
\ee
\end{small}
where $\hbar=t_1+t_2$. In the limit $t_1=t_2=0$, this operator gives the ordinary first Chern class (\ref{fcc}).
The last formula for $c_{1}({\cal{V}})$ identifies it with the well known hamiltonian of the quantum trigonometric Calogero-Moser system \cite{Ok1}. It  gives a first nontrivial element from the set of the cut-and-join operators studied for example in \cite{MMN1}-\cite{MMN2}.
The operators $c_{m}({\cal{V}})$ have even degrees and commute with each other, thus can be diagonalized  simultaneously. The characteristic classes in the equivariant cohomology always act diagonally on the classes of fixed points. Thus, the basis of $H^{\bullet}_{B}\Big(\coprod\limits_{n=0}^{\infty} \textrm{Hilb}_n \Big)$ consisting of the common eigenvectors for $c_{m}({\cal{V}})$ is given by the equivariant classes of the fixed points $\textrm{Hilb}_{n}^{B}$. The fixed points of $B$-action are naturally labeled by Young diagrams $\lambda$ with $|\lambda|=n$. The polynomial representing  the class $[\lambda]$ in the Fock space is the Jack polynomial with the parameter $\theta=-{t_1}/{t_2}$:
\be
j_{\lambda}(p_n)\in {\mathbb{C}}[p_1,p_2,...]\otimes {\mathbb{C}}[t_1,t_2]
\ee
The eigenvalues of $c_{1}$ for $j_{\lambda}(p_n)$ is the character of the fiber over the fixed point $\lambda$:
\be
c_{1}({\cal{V}})\,j_{\lambda}(p_n)= \sum\limits_{(i,j)\in \lambda}\, \Big( (i-1) t_1+(j-1) t_2 \,\Big)  j_{\lambda}(p_n)
\ee
The study of the Hilbert schemes gets simpler if one considers the embedding to a bigger space:
\be
\textrm{Hilb}_n \hookrightarrow {\cal{M}}(n,r)
\ee
where ${\cal{M}}(n,r)$ is the moduli space of framed, rank $r$ torsion free sheafs  on ${\mathbb{P}}^2$ with fixed Chern class $c_2=n$. There is a well known isomorphism between this space and the moduli space of the instantons on ${\mathbb{P}}^2$ and to abbreviate the notation, we refer to ${\cal{M}}(n,r)$ as the instanton moduli space in the text.  The framing of sheaf ${\cal{S}}$, by definition, is the choice of the isomorphism:
\be
\left.{\cal{S}}\right|_{\mathbb{P}}\rightarrow {\cal{O}}^{r}_{\mathbb{P}}
\ee
where $\mathbb{P}$ is a fixed line, usually considered as a line at infinity of a plane ${\mathbb{C}}^2\subset {\mathbb{P}}^2$. The moduli space carries the natural action of the group:
\be
\label{torss}
G=GL(r)\times GL(2)
\ee
where the second factor acts on ${\mathbb{C}}^2$ and the first one changes the framing. We denote by $A$ and $B$ the maximal tori of the first and the second factors in (\ref{torss}) respectively, and denote $T=A\times B $ the torus of the whole group.
%Note that $A\simeq ({\mathbb{C}}^{\ast})^2$ is the same group described above, and its action on the plane is given by (\ref{pa}).
Following \cite{OM} in section~\ref{Rmsec} we define a canonical map
\be
{\cal{R}} : H_{T}^{\bullet}\Big( \coprod\limits_{n=0}^{\infty} {\cal{M}}^{A}(n,r)\Big)\otimes \mathbb{C}(\frak{t}) \longrightarrow H_{T}^{\bullet}\Big(\coprod\limits_{n=0}^{\infty} {\cal{M}}^{A}(n,r)\Big) \otimes \mathbb{C}(\frak{t})
\ee
(Here the tensoring with rational functions  $\mathbb{C}(\frak{t})$ on the Lie algebra $\frak{t}$ of $T$ is the localization.)
For a fixed set we have:
\be
\label{fprt}
{\cal{M}}^{A}(n,r)=\coprod\limits_{n_1+...+n_r=n} \,\textrm{Hilb}_{n_1}\times...\times\textrm{Hilb}_{n_r}
\ee
The action of the first and the second component of (\ref{torss}) on the moduli space commute, thus $T/A=B$ acts on (\ref{fprt}). We have
$$
H_{T}^{\bullet}\Big(\coprod\limits_{n_1+...+n_r=n} \,\textrm{Hilb}_{n_1}\times...\times\textrm{Hilb}_{n_r}\Big) = {\cal{F}}^{\,\otimes r}\otimes \mathbb{C}[t_1,t_2,u_1,...,u_r]
$$
and, the map ${\cal{R}}$ can be considered as an operator:
\be
{\cal{R}}(u_1,...,u_r):\,{\cal{F}}^{\, \otimes r}\longrightarrow{\cal{F}}^{\, \otimes r}
\ee
acting in the tensor powers of the Fock space and depending on the equivariant characters $u_i$, $i=1...r$ of the torus $A$ as parameters.
The case that will be important to us is $r=2$, when we have the operator:
\be
\label{R}
{\cal{R}}(u): {\cal{F}}\otimes {\cal{F}}\longrightarrow{\cal{F}}\otimes {\cal{F}}
\ee
with $u=u_1-u_2$. The operators corresponding to $r>2$ are certain products of (\ref{R}).

By the construction,  ${\cal{R}}(u)$
satisfies the following Quantum Yang-Baxter Equation (QYBE) in ${\cal{F}}^{\otimes 3}$:
\be
{\cal{R}}_{12}(u) {\cal{R}}_{13}( u+ v) {\cal{R}}_{23}(v)={\cal{R}}_{23}(v) {\cal{R}}_{13}( u+v) {\cal{R}}_{12}(u)
\ee
where the operator ${\cal{R}}_{12}(u)$ acts as ${\cal{R}}(u)$ in the first two tensor components of ${\cal{F}}^{\otimes 3}$ and as identity in the last one, i.e. ${\cal{R}}_{12}(u)={\cal{R}}(u)\otimes 1$ and similarly for the other indices.

Given a solution of QYBE, we can apply the whole power of the quantum inverse scattering method \cite{QISMfirst}-\cite{QISMlast}, to identify  ${\cal{F}}$ with the Hilbert space of certain quantum integrable system. The set of commuting quantum hamiltonians for this system can be extracted from $R$-matrix as its specific matrix elements \cite{Sklyanin}. For example, in our case, we consider the following subspace in ${\cal{F}}^{\otimes 2}$:
$$
{\cal{F}}_{vac}=\textsf{vac}\otimes{\cal{F}}
$$
Here $\textsf{vac}$ is the \textit{vacuum vector} in ${\cal{F}}$, corresponding to $1$ in ${\mathbb{C}}[p_1,p_2,...]\otimes {\mathbb{C}}[t_1,t_2]$.
Let ${\cal{T}}(u):{\cal{F}}_{vac}\rightarrow{\cal{F}}_{vac}$ be the matrix element of ${\cal{R}}(u)$ for this subspace. Obviously, ${\cal{F}}_{vac}$ is isomorphic to a copy of the Fock space, such that ${\cal{T}}(u)$ is some operator in ${\mathbb{C}}[p_1,p_2,...]\otimes {\mathbb{C}}[t_1,t_2]$.
As direct consequence of QYBE \cite{Sklyanin}, we obtain
\be
\label{qis}
[{\cal{T}}(u),{\cal{T}}(v)]=0
\ee
We normalize $R$-matrix such that near $u=\infty$ the Taylor expansion has the form:
\be
\label{ex1}
{\cal{T}}(u)=1+\sum\limits_{n=1}^{\infty}\,{\cal{H}}_{n} u^{-n}
\ee
Then (\ref{qis}) implies that the coefficients  commute $[{\cal{H}}_{n}, {\cal{H}}_{m}]=0$ and thus can be identified with a set of commuting integrals of motion of some integrable system. Our case, in particular, corresponds to the trigonometric Calogero-Moser system with infinite number of particles, as for example, shown in~\cite{Ok1}.

The explicit form of the $R$-matrix, provides the formulae for the coefficients ${\cal{H}}_{n}$ directly in the form of boson operators. For example, the first coefficients take the form:
\be
\nonumber
\begin{array}{l}
{\cal{H}}_{1}=-\hbar \sum\limits_{n=1}^{\infty} \alpha_{-n}\alpha_{n},\\
\\
{\cal{H}}_{2}=\dfrac{\hbar}{2} \sum\limits_{m,n=1}^{\infty}\Big(  t_1 t_2 \alpha_{-m} \alpha_{-n} \alpha_{n+m}-\alpha_{-m-n} \alpha_{m} \alpha_n\Big)+\\
\\
+\dfrac{\hbar^2}{2} \Big(\sum\limits_{n=1}^{\infty} \alpha_{-n}\alpha_n\Big)^{2} +\dfrac{\hbar^2}{2} \sum\limits_{n=1}^{\infty} n \alpha_{-n}\alpha_n
\end{array}
\ee

 On the other hand, as follows from the definition of the $R$-matrix (see section \ref{rexpans}), ${\cal{T}}(u)$ must be diagonal in the basis of fixed points $j_{\lambda}$. The corresponding eigenvalues have the form:
\be
{\cal{T}}(u)\,j_{\lambda} = \left.\dfrac{ e( {\cal{V}}\otimes u ) }{e({\cal{V}}\otimes u \otimes \hbar)}\right|_{\lambda} j_{\lambda}
\ee
where $e(L)\in H_{T}^{\bullet}(X)$ stands for the $T$-equivariant  Euler class of the bundle $L$, and we denoted by the characters $\hbar=t_1+t_2$ and $u$ the trivial equivariant line bundles with the corresponding scaling action of the torus $T$ on their fibers. The tautological bundle ${\cal{V}}$ will be identified with certain subbundle of the normal bundle to $\textrm{Hilb}_{n}$ in ${\cal{M}}(n,r)$.

By virtual splitting, assume that the tautological bundle ${\cal{V}}$ splits to the sum of line bundles. Let us enote by $x_i$, $i=1...\textrm{rk} {\cal{V}}$ the corresponding Chern roots. Then, for the equivariant Euler class we obtain:
\be
\label{eul}
\dfrac{ e(  {\cal{V}}\otimes u  ) }{e({\cal{V}}\otimes u\otimes \hbar)}=\dfrac{\prod\limits_{i=1}^{\textrm{rk }{\cal{V}} } (u+x_i)}{\prod\limits_{i=1}^{\textrm{rk}  {\cal{V}} } (u+\hbar+x_i)}
\ee
The Chern classes of ${\cal{V}}$ are the elementary symmetric polynomials in the Chern roots, i.e. the generating function is:
\be
\label{cr}
\prod\limits_{k=1}^{\textrm{rk}  {\cal{V}}} (1+z x_{i})=\sum\limits_{i=1}^{\textrm{rk}  {\cal{V}}} \, z^i c_{i}({\cal{V}})
\ee
Therefore expanding (\ref{eul}) in power series  at $u=\infty$, using $(\ref{cr})$ we obtain:
\be
\label{ex2}
{\cal{T}}(u)=1-\dfrac{\hbar\, \textrm{rk}  {\cal{V}} }{u}+\dfrac{\hbar c_1 ({\cal{V}}) + \hbar^2 \textrm{rk}  {\cal{V}}(\textrm{rk}  {\cal{V}}+1)/2}{u^2} +O(u^{-2})
\ee
so that in general the coefficient of $u^{-n}$  in the expansion is some universal polynomial in Chern classes $c_{i}({\cal{V}})$  and rank $\textrm{rk}  {\cal{V}}$. Now, comparing the coefficients of expansion in (\ref{ex2}) with same expansion (\ref{ex1}) give us the explicit formulas for operators $c_n ({\cal{V}})$ in terms of the Nakajima boson operators. For example, we obtain:
\be
\textrm{ rk} {\cal{V}} =\sum\limits_{k=1}^{\infty} \alpha_{-k} \alpha_k
\ee
Analogously, the coefficient in $u^{-2}$ gives Lehn's formula (\ref{lf}).

This article has the following structure: section \ref{sec1} contains a short introduction to the cohomology of Hilbert schemes and  instanton moduli. Following \cite{OM} (section 3) in  section \ref{sten}  we define the stable map in the cohomology of symplectic  $T$-space. We describe the tautological bundle on the Hilbert schemes as the component of the normal bundle in the instanton moduli space. In section \ref{Rmsec} we define the $R$-matrix and compute it explicitly in the case of cotangent bundles to flag varieties. We describe the factorization procedure for the $R$-matrices. This procedure gives us certain product formulae for the $R$-matrices associated with $A_n$ -quiver varieties and instanton moduli.  In section we compute explicitly unknown factors in the product formulas of section \ref{Rmsec}, using the fusion procedure from QISM.  This finishes the computations of R-matrices in the case of $A_n$ varieties and instanton moduli spaces. In the last section \ref{expsec} we study the matrix elements and the power series expansion of instanton $R$-matrix. In particular we rederive the Lehn's formula for the first Chern class.

\vspace{2mm}

\noindent
\textbf{Acknowledgements:} I would like to thank A. Okounkov for suggesting this project to me, countless explanations, discussions of the subject and his interest to this work. I am grateful to Michael McBreen and Andrei Negut for valuable discussions, and all other participants of Columbia mathematical physics seminar. This work was supported in part by RFBR grants 12-01-33071 mol-a-ved, 12-02-00594 and  12-01-00482.
%%%%%%%%%%%%%%%%%%%%%%%%%%%%%%%%%%%%%%%%%%%%%%%%%%%%%%%%%%%%%%%%%%%%%%%%%%%%%%%%%%%%%%%%%%
%%%%%%%%%%%%%%%%%%%%%%%%%%%%%%%%%%%%%%%%%%%%%%%%%%%%%%%%%%%%%%%%%%%%%%%%%%%%%%%%%%%%%%%%%%
\section{Hilbert Schemes \label{sec1}}
\subsection{Equivariant cohomology of Hilbert Schemes}
The Hilbert scheme $\textrm{Hilb}_n$ of $n$ points in the plane ${\mathbb{C}}^2$ is the configuration space of $n$-tuples of points on  ${\mathbb{C}}^2$. It can be described as the space of polynomial ideals ${\cal{J}}\subset {\mathbb{C}}[x,y]$ of complex codimension $n$:
$$
\dim_{{\mathbb{C}}} \, {\mathbb{C}}[x,y]/{\cal{J}}=n
$$
The Hilbert scheme $\textrm{Hilb}_n$ has rich, well understood geometry \cite{Lehn,Nak1}. It is a nonsingular, irreducible, quasiprojective algebraic variety of dimension $2n$. The symmetries of ${\mathbb{C}}^2$ lift to the Hilbert scheme. In particular, the scaling of coordinates in ${\mathbb{C}}^2$:
$$
(z_1,z_2) \cdot(x,y) =(z_1 x, z_2 y)
$$
induces the action of algebraic torus $B={({\mathbb{C}}^{\ast})}^2$ on $\textrm{Hilb}_n$. The $B$-equivariant cohomology of the Hilbert scheme can be described in the framework of the boson Fock space formalism. Let ${{\frak h}}$ be the Heisenberg algebra generated by bosons $\alpha_n$ satisfying the following commutation relations:
$$
[\alpha_n, \alpha_m ]=n \delta_{m+n}
$$
The Fock space ${\cal{F}}$ is an infinite dimensional irreducible representation of the Heisenberg algebra in the space of symmetric polynomials on infinite number of variables:
$$
{\cal{F}}={\mathbb{C}}[p_1,p_2,p_3,...]
$$
with the following action of the bosons:
\be
\alpha_{-n}=p_n, \ \ \ \alpha_{n}= n \dfrac{\partial}{\partial p_n}, \ \ \textrm{for} \ \ \ n>0
\ee
The element $\textsf{vac} \in {\cal{F}}$ corresponding to $1\in {\mathbb{C}}[p_1,p_2,p_3,...]$ is called vacuum vector. The operators $\alpha_{-n}$ for $n>0$ are called creation operators, as they create the state $\alpha_{-n} \textsf{vac} =p_n$ from the vacuum. The operators $\alpha_{n}$ are called the annihilation operators. They kill the vacuum vector $\alpha_n \textsf{vac}=0$.

The Fock space has a natural basis indexed by the partitions $\lambda$:
\be
\label{vec}
| \lambda \rangle= \dfrac{1}{{\frak{d}}(\lambda)} \prod\limits_{i=1}^{l(\lambda)}{\alpha_{-\lambda_i}}\, \phi= \dfrac{1}{{\frak{d}}(\lambda)} \prod\limits_{i=1}^{l(\lambda)}{p_{\lambda_i}}\,
\ee
where $l(\lambda)$ is the length of the partition, and
$$
{\frak{d}}(\lambda)=|\textrm{Aut}(\lambda)| \prod\limits_{i}^{l(\lambda)} \lambda_i
$$
is the standard combinatorial factor.

The action of $B$ on the Hilbert scheme is hamiltonian and thus the equivariant cohomology is a tensor product \cite{Brion}:
\be
\label{coh}
{\cal{H}}=H^{\bullet}_{B}\Big(\coprod\limits_{n=0}^{\infty} \, \textrm{Hilb}_n \Big)=\bigoplus\limits_{n=0}^{\infty} H^{\bullet}_{B}(\textrm{Hilb}_n) ={\cal{F}}\otimes_{{\mathbb{C}}} {\mathbb{C}}[t_1,t_2]
\ee
where $t_1$ and $t_2$ are the equivariant parameters corresponding to the weights of $B$. The whole space (\ref{coh}) is a free module over
the equivariant cohomology of a point
$$
{\mathbb{C}}[t_1,t_2]=H^{\bullet}_{B}(\textrm{pt})
$$
The usual cohomology ring is a quotient \cite{Brion}:
\be
H^{\bullet} \Big(\coprod\limits_{n=0}^{\infty} \, \textrm{Hilb}_n \Big)=H^{\bullet}_{B}\Big(\coprod\limits_{n=0}^{\infty} \, \textrm{Hilb}_n \Big)/(t_1,t_2)\simeq {\cal{F}}
\ee
Thus, all  formulas for the ordinary cohomology can be obtained from the equivariant theory in the limit $t_1=t_2=0$.

The vacuum vector $\textsf{vac}$ corresponds to the unit in $H^{\bullet}(\textrm{Hilb}_0)\simeq {\mathbb{C}}$. The subspace in ${\cal{F}}$ corresponding to $H^{\bullet}(\textrm{Hilb}_n)$ is spanned by the vectors (\ref{vec}) with $|\lambda|=n$.
For example, we have:
$$
H^{\bullet}(\textrm{Hilb}_4)=\textrm{Span}\Big( p_{1}^4,\, p_{2}p_{1}^2, \, p_{3}p_{1},\, p_{2}^2, p_{4} \Big)
$$
The cohomological degree of the monomial is $\textrm{deg}(p_{n})=2(n-1)$. For example, for $\textrm{Hilb}_4$ we obtain:
$$
\begin{array}{c}
H^{0}(\textrm{Hilb}_4)=\textrm{Span}( p_{1}^4 ),\ \ H^{2}(\textrm{Hilb}_4)=\textrm{Span}( p_{2} p_{1}^2 ), \\
\\
 H^{4}(\textrm{Hilb}_4)=\textrm{Span}( p_{3} p_{1},\, p_2^2  ), \ \ H^{6}(\textrm{Hilb}_4)=\textrm{Span}( p_4 )
\end{array}
$$
Geometrically, the vector $|\lambda \rangle$ in (\ref{vec}) corresponds to the class of the subvariety of $\textrm{Hilb}_{|\lambda|}$ given by the union of schemes of length $\lambda_1,\,\lambda_2,\, ..., \lambda_{l(\lambda)}$ supported at $l(\lambda)$ distinct points in the plane ${\mathbb{C}}^2$.

The ring structure of ${\cal{H}}$ is more complicated. Of course, by obvious dimensional reasons it does not coincide with the ring structure in ${\mathbb{C}}[p_1,p_2,p_3,...]$. For every element $\gamma \in {\cal{H}}$ one can define the operator of the cup product:
$$
m(\gamma) \in \textrm{End} \,( {\cal{H}}), \ \ \ m(\gamma)\cdot \beta \mapsto \gamma \cup \beta
$$
For example, the vector $|1^n\rangle=|{1,1,...,1}\rangle =p_1^n$, having degree zero, corresponds to the identity operator in ${\cal{H}}$ :
\be
\left.m\Big(\,| 1^n \rangle\,\Big)\right|_{H^{\bullet}_{B}(\textrm{Hilb}_n)} =\textrm{Id}
\ee
The action of the operators $m(\gamma)$ on ${\cal{H}}$ can be described by certain infinite sums of bosons $\alpha_n$ known as generalized cut-and join-operators. The algebra of these operators and the correspondence $\gamma \rightarrow m(\gamma)$ is well studied for example in \cite{MMN1,MMN2}.

As a ring ${\cal{H}}$ is generated by the equivariant Chern classes of the tautological bundle:
\be
{\cal{V}}={\cal{O}}/{\cal{J}}\rightarrow \coprod\limits_{n=0}^{\infty} \textrm{Hilb}_{n}
\ee
with the fiber ${\mathbb{C}}[x,y]/{\cal{J}}$ over the point $[\cal{J}]$, such that $\left.{\cal{V}}\right|_{\textrm{Hilb}_n}$ is the rank $n$ bundle over $\textrm{Hilb}_n$. The first Chern class of this bundle is represented by the unique element of degree 2 in $H^{2}(\textrm{Hilb}_n)$:
\be
c_{1}({\cal{V}})=|2,1^{n-2}\rangle =p_2 p_1^{n-1} \in H^{2}_{B}(\textrm{Hilb}_n)
\ee
The cup product with this class is described by the following cut-and-join operator:
$$
\begin{array}{r}
c_{1}({\cal{V}})=\dfrac{1}{2}\sum\limits_{m,n=1}^{\infty}\Big( t_1 t_2 \alpha_{-m} \alpha_{-n} \alpha_{n+m} - \alpha_{-m-n} \alpha_{n} \alpha_{m}\Big)+
\dfrac{\hbar}{2}\, \sum\limits_{n=1}^{\infty}\,( n-1) \alpha_{-n}\alpha_{n}
\end{array}
$$
This operator coincides with the well known Hamiltonian of the trigonometric Calogero-Moser-Sutherland system \cite{Ok1}.
\subsection{Embedding to the instanton moduli space \label{embd}}
%\verb"Here should be explained the ADHM description"
Let ${\cal{M}}(n,r)$ be the moduli space of framed, torsion free sheaves ${\cal{S}}$ on ${\mathbb{P}}^2$ with rank $r$ and $c_2({\cal{S}})=n$. By a framing of a sheaf ${\cal{S}}$
we mean the choice of an isomorphism:
$$
\left.{\cal{S}}\right|_{{\mathbb{P}}}\rightarrow {\cal{O}}^{\oplus r}_{L_{\infty}}
$$
where ${L_{\infty}}\subset {{\mathbb{P}}}^2$ is a fixed line, usually considered as the infinity line of ${\mathbb{C}}^2\subset {\mathbb{P}}^2$.
The moduli space ${\cal{M}}(n,r)$ comes with the natural action of a group:
\be
\label{gr}
G=\textrm{GL}(r)\times \textrm{GL}(2)
\ee
The first factor acts by changing the framing and the second one comes from the action on ${\mathbb{C}}^2$.

The space ${\cal{M}}(n,r)$ can be described as a Nakajima variety corresponding to a quiver with one vertex and one loop \cite{Nak1}. This description is, of course, the well known ADHM construction of instantons \cite{ADHM}.
Let $V={\mathbb{C}}^n$ and $W={\mathbb{C}}^r$ be spaces with fixed dimensions corresponding to the vertices of the quiver, then
\begin{small}
\be
\label{adhm}
{\cal{M}}(n,r)=\left.\left\{ (B_1,B_2,i,j) \left| \begin{array}{l} [B_1,B_2] + i j=0\\  \textrm{there is no subspace} \ \  S\subset V, \textrm{such that}\ \ \\ B_{\alpha}(S)\subset S, \,  \alpha=1,2 \ \ \textrm{and} \ \  \textrm{Im}(i)\subset S  \end{array}\right.  \right\}\right/\textrm{GL}(V)
\ee
\end{small}
where $B_1, B_2 \in \textrm{End}(V)$, $i \in \textrm{Hom}(W,V)$, $j\in \textrm{Hom}(V,W)$, with the action of $g\in \textrm{GL}(V)$ given by:
$$
g \cdot (B_1,B_2,i,j)=(g B_1 g^{-1},g B_2g^{-1} ,g i,j g^{-1})
$$
For $r=1$ the stability condition gives $j=0$, and we arrive to the standard ADHM description of the Hilbert schemes of points on a plane \cite{Nak1}:
\be
{\cal{M}}(n,1) = \textrm{Hilb}_{n}
\ee

The moduli space is symplectic. The symplectic form induced by the skew-symmetric form on the representation space of the quiver:
\be
\label{sf}
\omega\Big( (B_1,B_2,i,j), (B_1^{\prime},B_2^{\prime},i^{\prime},j^{\prime}) \Big)=\textrm{tr}\Big(B_1 B_2^{\prime}-B_2 B_1^{\prime}+i j^{\prime}-i^{\prime} j\Big)
\ee
The action of the first factor of $(\ref{gr})$ on the moduli space is given by
\be
\label{gwact}
g_w\cdot (B_1, B_2, i,j)=(B_1,B_2, i g_w^{-1}, g_w j),
\ee
with $g_w\in  \textrm{GL}(r)$. The second factor $ \textrm{GL}(2)$ acts by
\be
g_l \cdot (B_1,B_2)= (a B_1+b B_2, c B_1+d B_2),
\ee
with
$$
\label{el}
g_l=\left(\begin{array}{cc} a& b\\ c & d \end{array}\right) \in \textrm{GL}(2)
$$
The action of the second factor $ \textrm{GL}(r)$, obviously, preserves the symplectic form (\ref{sf}).  The first factor $ \textrm{GL}(2)$ scales it i.e. the subgroup $ \textrm{SL}(2)\subset \textrm{GL}(2)$ is the stabilizer of $\omega$.

Let $A$ and $B$ are the maximal tori of the first and the second factor of (\ref{gr}) respectively, such that the full torus acting on the moduli space is:
\be
\label{tors}
T=A\times B
\ee
The torus $B\simeq {\mathbb{C}^{\ast}}^2$ and we denote by $t_1$, $t_2$ the generators of the character group $c(B)\simeq {\mathbb{Z}}^2$. As we mentioned above, the action of $B$ scales the symplectic form, thus ${\mathbb{C}}\, \omega$ is a one dimensional representation of $B$ with character $\hbar=t_1+t_2$.

The general fact about the quiver varieties is that, the fixed points set of torus actions preserving the symplectic form are also quiver varieties.
The situation is different, however, for the vertex part $A$ and the loop part $B$ of the torus (\ref{tors}). The fixed points of the $A$-action on is the union of smaller quiver varieties of the same type:
$$
{\cal{M}}^{A}(n,r)=\coprod\limits_{n_1+...+n_r=n}\, {\cal{M}}(n_1,1)\times...\times {\cal{M}}(n_r,1)
$$
Thus, as was noted above:
\be
\label{emb}
{\cal{M}}^{A}(n,r)=\coprod\limits_{n_1+...+n_r=n}\,\textrm{ Hilb}_{n_1}\times...\times \textrm{Hilb}_{n_r}
\ee
The sheaf  ${\cal{S}}\in {\cal{M}}^{A}(n,r)$ corresponds to the sum of the ideal sheaves:
$$
{\cal{S}}=I_{1}\oplus...\oplus I_{r}
$$
In particular, we have:
\be
H^{\bullet}_{T}\Big(  \coprod\limits_{n=0}^{\infty} {\cal{M}}^{A}(n,r)\Big)\simeq{\cal{F}}^{\otimes r}\otimes \mathbb{C}[t_1,t_2,u_1,...,u_r]
\ee
where $u_i$ are the characters of $A$ and ${\mathbb{C}}[t_1,t_2,u_1,...,u_r]=H_{T}^{\bullet}(\textrm{pt})$.
The decomposition (\ref{emb}) gives an embedding of $\textrm{Hilb}_{n}$ to the moduli space as a component of the fixed set:
\be
\textrm{Hilb}_{n}=\textrm{Hilb}_{0}\times \textrm{Hilb}_{0}\times ... \times \textrm{Hilb}_{n} \subset {\cal{M}}^{A}(n,r)
\ee
This allows one to describe the topology of $\textrm{Hilb}_{n}$ in terms of this embedding. For example, in the next section we identify the tautological bundle on $\textrm{Hilb}_{n}$ with certain component of the normal bundle to $\textrm{Hilb}_{n}$ in ${\cal{M}}(n,r)$.

Let us consider a general subtorus corresponding of the loop part ${\mathbb{C}}^{\ast}\subset B$ defined as the elements (\ref{el}) of the form:
$$
g_l=\left(\begin{array}{cc} z& 0\\ 0 & z^{-1} \end{array}\right)
$$
The action of this subgroup gives the decomposition of $V$ and $W$ to weight spaces:
\be
V=\bigoplus\limits_{k\in {\mathbb{Z}}} V^{(k)}, \ \ \ W=\bigoplus\limits_{k\in {\mathbb{Z}}} W^{(k)}
\ee
If $(B_1,B_2,i,j)$ represents a class of the fixed point for this torus, then all elements must have the following block form:
$$
B_{i}=\bigoplus\limits_{k\in {\mathbb{Z}}} B_{i}^{(k)}, \ \ i=\bigoplus\limits_{k\in {\mathbb{Z}}} i^{(k)}, \ \ j=\bigoplus\limits_{k\in {\mathbb{Z}}} j^{(k)}
$$
with:
\be
\nonumber
B_{1}^{(k)}\in \textrm{Hom}(V^{(k)},V^{(k+1)}), \ \ B_{2}^{(k)}\in \textrm{Hom}(V^{(k+1)},V^{(k)}),\\
\nonumber
\\ \nonumber
 i^{(k)}\in \textrm{Hom}(W^{(k)},V^{(k)}), \ \  j^{(k)}\in \textrm{Hom}(V^{(k)},W^{(k)}) \ \ \
\ee
Note that these data defines the representation of $A_{\infty}$ quiver, and the stability condition in (\ref{adhm}) translates to the stability condition for  $A_{\infty}$ quiver variety. Thus, for a set of fixed points we obtain:
\be
{\cal{M}}^{B}(n,r)=\coprod\limits_{|v|=n, |w|=r} {\cal{A}}_{\infty}(v,w)
\ee
where ${\cal{A}}_{\infty}(v,w)$ is $A_{\infty}$-quiver variety with dimensions $\dim V^{(k)}=v_k$, $\dim W^{(k)}=w_k$ satisfying $|v|=\sum v_k=n$, and $|w|=\sum w_k=r$.
\section{Stable envelope \label{sten}}
In this section, following \cite{OM}, we recall the definition of \textit{the stable map} playing important role in the description of $R$-matrix. Here we consider this construction in general. Concrete examples for the flag manifolds and the instanton moduli spaces will be considered in the next section.

Assume that a pair of algebraic tori  $A\subset T$ acts on the symplectic, quasiprojective algebraic variety $X$. This action induces the action on $H^{0}(\Omega^2_{X})$. Assume that the induced action of $T$ on $H^{2}(X)$ scales the symplectic form $\omega$. It implies that the one-dimensional
subspace ${\mathbb{C}}\omega\subset H^{2}(X)$ is a subrepresentation of $T$. We denote by $\hbar$ its character. Assume, that the action of the smaller torus $A$ preserves $\omega$.

%We assume that $X$ has a symplectic form $\omega\in H^{0}(\Omega^2_{X})$, which is an eigenvector of the induced $T$-action, and $\omega$ is fixed by $A$.
%Let $i^{\ast} : H^{\bullet}_{T}(X)\rightarrow H^{\bullet}_{T}(X^{A})$ is the map induced by the inclusion $i:=X^{A}\rightarrow X$ of the fixe .
Our goal in this section is to describe the natural map defined in \cite{OM}:
$$
\textrm{Stab}_{\fC}: H^{\bullet}_{T}(X^{A}) \rightarrow H^{\bullet}_{T}(X)
$$
depending on chamber $\fC$ in the Lie algebra $a_{\mathbb{R}}$.  For a fixed cycle $\gamma \in H^{\bullet}_{T}(X^{A})$ the element
$\textrm{Stab}_{\fC}(\gamma)\in H^{\bullet}_{T}(X)$ is called \textit{stable envelope} of $\gamma$.

%%%%%%%%%%%%%%%%%%%%%%%%%%%%%%%%%%%%%%%%%%%%%%%%%%%%%%%%%%%%%%%%%%%%%%%%%%%%%%%
%%%%%%%%%%%%%%%%%%%%%%%%%%%%%%%%%%%%%%%%%%%%%%%%%%%%%%%%%%%%%%%%%%%%%%%%%%%%%%%

\subsection{Chamber decomposition \label{cdec}}

Let $A\simeq ({\mathbb{C}}^{\ast})^r $ be an algebraic torus of rank $r$. Let
\be
\label{coc}
 c(A)=\{\, A \rightarrow {\mathbb{C}}^{\ast} \,\}\simeq{{\mathbb{Z}}}^r, \ \ \ t(A)=\{\, {\mathbb{C}}^{\ast} \rightarrow A \,\}\simeq{{\mathbb{Z}}}^r
\ee
be the group of characters and cocharacters respectively.  We define the real part of the Lie algebra and its dual as:
\be
\label{rlie}
a_{{{\mathbb{R}}}}=t(A)\otimes_{{\mathbb{Z}}}{\mathbb{R}}\simeq{\mathbb{R}}^r\subset \textrm{Lie}(A), \ \ \ a_{{{\mathbb{R}}}}^{\ast}=c(A)\otimes_{{\mathbb{Z}}}{\mathbb{R}}\simeq{\mathbb{R}}^r
\ee
The natural pairing $t(A)\times c(A)\rightarrow {{\mathbb{Z}}}$ linearly extends to the pairing for the real Lie algebra:
$$\langle \ \ , \ \ \rangle : a_{{{\mathbb{R}}}} \times a_{{{\mathbb{R}}}}^{\ast} \rightarrow {{\mathbb{R}}}$$

\noindent
\textbf{Definition:}
Let $X^{A}$  be the fixed set of $A$. The normal bundle $N$ to $X^A$ in $X$ has a natural structure of an $A$-module and  splits to the direct sum of complex, one-dimensional, irreducible components. The subset $\Delta\subset a_{\mathbb{R}}^{\ast}$ consisting of the characters appearing in $N$ is called \textit{root system } of $A$.

A weight $\alpha \in a_{{{\mathbb{R}}}}^{\ast}$ defines a hyperplane in $a_{{{\mathbb{R}}}}$:
\be
\textrm{ker}_\alpha=\{v\in a_{{{\mathbb{R}}}}: \langle \alpha, v \rangle =0\}
\ee
The  hyperplanes corresponding to the roots partition $a_{{{\mathbb{R}}}}$ into the set of open chambers:
\be
a_{{\mathbb{R}}} \setminus \bigcup\limits_{\alpha \in \Delta}\, \textrm{ker}_{\alpha }=\coprod\limits_{i} \fC_{i}
\ee

The  hyperplanes $\textrm{ker}_{\alpha}$, clearly,  are walls of the chambers. In general, the walls $\textrm{ker}_{\alpha}$ define a stratification of the space $a_{{\mathbb{R}}}$ by the chain of sets: the set of points that do not lie on any wall (these are chambers), the set of points lying on exactly one wall, the set of points lying on the intersection of two walls and so on.

The stratification of $a_{{\mathbb{R}}}$ encodes the information about the $A$-action on $X$. Indeed, consider a cocharacter $\sigma: {\mathbb{C}}^{\ast} \rightarrow A$. It defines certain ${\mathbb{C}}^{\ast}$-action on $X$. If $\sigma$ does not belong to some wall i.e. is inside one of the chambers then,
${\mathbb{C}}^{\ast}$ - action has the same set of the fixed points $X^{{\mathbb{C}}^{\ast}} = X^{A}$. Assume now, that $\sigma$ lies on exactly one wall $\ker_{\alpha}$. Then torus ${\mathbb{C}}^{\ast}$ acts trivially on the component of the normal bundle corresponding to the character $\alpha$. Therefore,  ${\mathbb{C}}^{\ast}$ preserves the corresponding direction in $X$ and the fixed set $X^{{\mathbb{C}}^{\ast}}$ gets larger then $X^{A}$.   Thus, the stratification of $a_{{{\mathbb{R}}}}$ by the walls $\ker_{\alpha}$ corresponds to the types of the fixed sets $X^{{\mathbb{C}}^{\ast}}$ arising from different choice of the cocharacters  $\sigma: {\mathbb{C}}^{\ast} \rightarrow A$. In the extreme case $\sigma=0$, corresponding to the intersection of all walls, we have $X^{{\mathbb{C}}^{\ast}}=X$.
%the roots are again of the form $a_{i}-a_{j}$ and correspond to the standard roots of $GL(r)$, and the chambers are the Weyl chambers.  $\fC$
\subsection{Stable leaves and slopes}
\noindent
\textbf{Definition:}
Let us fix some chamber $\fC \subset a_{\mathbb{R}}$ and let $\sigma\in \fC$ be a cocharacter. We say that the point $x \in X$ is $\fC$\textit{-stable} if the following limit exists:
\be
\lim_{\fC} x\stackrel{\textrm{def}}{=} \lim\limits_{z\rightarrow 0} \sigma(z)\cdot x \in X^{A}
\ee
This definition does not depend on the choice of $\sigma$ in the chamber $\fC$, which explains the notation $\lim\limits_{\fC}$.
For a component $Z$ of the fixed set $ X^{A}$ we define its  \textit{stable leaf} as the set of stable points "attracting" to $Z$:

\noindent
\textbf{Definition:}
\be
\textrm{Leaf}_{\fC}(Z)=\{ x|\lim_{\fC} x  \in Z \}
\ee
The choice of a chamber $\fC$ defines a partial order on the components $Z\subset X^{A}$. We say that:
 $$Z_1 \succeq Z_{2} \ \  \Leftrightarrow \ \ \overline{\textrm{Leaf}_{\fC}(Z_1)} \cap Z_2\neq \emptyset.$$ Using this ordering we define the  \textit{stable slope} of a component  $Z\subset X^{A}$ as follows:

\noindent
\textbf{Definition:}
\be
\textrm{Slope}_{\fC}(Z)=\coprod\limits_{Z^{\prime} \preceq\, Z} \textrm{Leaf}_{\fC}(Z^{\prime}).
\ee
\noindent

\subsection{Polarization}
 Let $Z\subset X^{A}$ be a component of the fixed set. The choice of a chamber~$\fC$ gives a decomposition of the normal bundle to $Z$ in $X$ into the weight spaces that are positive or negative on the chamber $\fC$:
$$
N_{Z}=N_{+}\oplus N_{-}
$$
Remind, that by our assumption the action of $A$ preserves the symplectic form $\omega$ and $T$ scales it with character $\hbar$. Thus, we have:
\be
\label{polar}
(N_{+})^{\vee}=N_{-}\otimes \hbar
\ee
where for convenience we denoted by the same symbol $\hbar$  the trivial $T$-equivariant line bundle over $Z$ with the action of $T$ on its fiber corresponding to the character $\hbar$.

Assume that $\alpha_{i}$, $i=1...\textrm{codim}(Z)/2$ are the weights of $N_{+}$ then,  the $A$-weights of $N_{-}$ are given by $(-\alpha_{i})$. Therefore, the $A$-equivariant Euler class of $N_{Z}$ (with a sign) is a perfect square:
\be
\label{Ec}
\varepsilon^2=(-1)^{\textrm{codim}(Z)/2} e(N_Z)=\prod\limits_{i=1}^{\textrm{codim}(Z)/2} \alpha_{i}^2\ \
\ee

\noindent
\textbf{Definition:}
The \textit{polarization} of $Z$ is a formal choice of a sign in the square root of (\ref{Ec}):
\be
\left.\varepsilon\right|_{H^{\bullet}_{A}(\textrm{pt})}=\pm \prod\limits_{i=1}^{\textrm{codim}(Z)/2} \alpha_{i}
\ee
We say that the sign $\pm e(N_{-})\in {H^{\bullet}_{T}(Z)}$ is chosen according to the polarization if it restricts to $\varepsilon$ in $H^{\bullet}_{A}(Z)$.
\subsection{Tautological bundle \label{tb}}
Let $Z=\textrm{Hilb}_{n_1}\times \textrm{Hilb}_{n_2}$, $n_1+n_2=n$ be a component of the fixed set ${\cal{M}}^{A}(n,r)$ for $r=2$. Assume that some chamber $\fC$ in the Lie algebra of $A$ is fixed, then it defines a decomposition:
$$
N_Z=N_{-}\oplus N_{+}
$$
To describe these components explicitly, let us fix a cocharacter of $A$ such that ${\mathbb{C}}^{\ast}$ acts on ${\cal{M}}(n,r)$ by (\ref{gwact}) with a matrix:
$$
g_w=\left(\begin{array}{cc} z^{u_1}& 0\\ 0 & z^{u_2} \end{array}\right)
$$
Consider a point ${\cal{S}}={\cal{I}}\oplus {\cal{J}} \in Z $. The tangent space to the moduli space at this point is:
$$
T_{{\cal{S}}} {\cal{M}}(n,2)=\textrm{Ext}^{1}({\cal{S}},{\cal{S}})
$$
such that it is a sum of four terms:
\be
\label{ext}
T_{{\cal{S}}}\, {\cal{M}}(n,2)=\textrm{Ext}^{1}({\cal{I}},{\cal{I}})\oplus\textrm{Ext}^{1}({\cal{I}},{\cal{J}})
\oplus\textrm{Ext}^{1}({\cal{J}},{\cal{I}})\oplus\textrm{Ext}^{1}({\cal{J}},{\cal{J}})
\ee
The torus $A$ acts naturally on this spaces with characters $0$, $u_1-u_2$, $u_2-u_1$ and $0$ respectively. Fixing a chamber of the cocharacters taking negative values on $u=u_1-u_2$ we obtain the identification:
\be
N_{-}=\textrm{Ext}^{1}({\cal{I}}, {\cal{J}} )
\ee
Clearly, the spaces $\textrm{Ext}^{1}({\cal{I}},{\cal{I}})$, $\textrm{Ext}^{1}({\cal{J}},{\cal{J}})$ are the  tangent spaces to the components of $\textrm{Hilb}_{n_1}\times \textrm{Hilb}_{n_2}$. The remaining terms in (\ref{ext}) give the decomposition of the normal bundle to the fixed set into the positive and negative part:
$$
N_Z=N_-\oplus N_{+}=\textrm{Ext}^{1}({\cal{I}},{\cal{J}})\oplus \textrm{Ext}^{1}({\cal{J}},{\cal{I}})
$$
The following theorem describes the tautological bundle ${\cal{V}}$ in terms of the embedding  $\textrm{Hilb}_{n}\subset {\cal{M}}(n,2)$.

\noindent
\textbf{Theorem}
Let $Z=\textrm{Hilb}_{n}=\textrm{Hilb}_{0}\times \textrm{Hilb}_{n}\subset {\cal{M}}(n,2)$ be a component of the fixed set and  $N_{-}$ be the corresponding component of the normal bundle. We have the identification:
\be
N_{-}={\cal{V}}\otimes u
\ee
where the character $ u $ denotes the trivial line bundle over $\textrm{Hilb}_n$ with the action of $A$ on its fiber given by $u$.

\noindent
\textit{Proof.} In this case ${\cal{I}}\otimes {\cal{J}}\in \textrm{Hilb}_{0}\otimes \textrm{Hilb}_{n}$, therefore ${\cal{I}}={\cal{O}}$. We have
$$
N_-=\textrm{Ext}^{1}({\cal{O}}, {\cal{J}})=H^{1}({\cal{J}})
$$
Now, from the long exact sequence associated with the short one:
\be
0\rightarrow {\cal{J}} \rightarrow {\cal{O}} \rightarrow {\cal{O}} / {\cal{J}} \rightarrow 0
\ee
we obtain $H^{1}({\cal{J}})=H^{0}({\cal{O}} / {\cal{J}})={\cal{V}}$ as the cohomology of the structure sheaf ${\cal{O}}$ vanish.  As we discuss above, $A$ acts on $N_-$ with the character $u$ which gives the result of the theorem.

We will also need the explicit formula for the equivariant Euler class of $N_{-}$ given by the following fact:

\textbf{Theorem}
Let ${\cal{S}}=I_{\lambda}\oplus I_{\mu} \in {\cal{M}}^{A}(n,2)$ be a fixed point of $T$-action and $\lambda$, $\mu$ be the partitions labeling the fixed points one the Hilbert schemes.  Then we have the following explicit formula for the equivariant Euler class:
$$
\left.e(N_{-})\right|_{{\cal{S}}}=\prod\limits_{\Box \in\, \lambda} \Big(u+t_1\, l_\mu(\Box) -t_2 (a_{\lambda}(\Box)+1)  \Big)
\prod\limits_{\Box \in \,\mu} \Big( u-t_1\, (l_{\lambda}(\Box)+1) + t_2 \,a_{\mu}(\Box) \Big)
$$
where $a_\lambda(\Box)$ and $l_\lambda(\Box)$ are the arm and leg lengths of the box $\Box$ in the partition $\lambda$. If the box $\Box$ has the coordinates $i, j$ then they are defined as:
$$
a_\lambda(\Box)=\lambda_i-j, \ \ l_{\lambda}(\Box)=\lambda^{'}_{j}-i
$$
\noindent
The proof of this theorem is given, for example, in \cite{CarOk}.

Note, that in the case $Z=\textrm{Hilb}_{n}$, corresponding to $\lambda=\emptyset$, the above formula specializes to the Euler class of the tautological bundle over $\textrm{Hilb}_{n}$:
\be
\label{Nekf}
e(N_{-})=e({\cal{V}}\otimes u)=\prod\limits_{(i,j)\in \mu}\,\Big(u+t_1 (j-1)+t_2 (i-1)\Big)
\ee

\subsection{Stable envelope \label{se}}

The stable envelope is defined by the following theorem.

\textbf{Theorem} Under assumption above, there exists a unique map of $H_{T}^{\bullet}(\textrm{pt})$ modules:
$$
\textrm{Stab}_{\fC, \varepsilon} : H^{\bullet}_{T}(X^{A}) \rightarrow H_{T}^{\bullet}(X)
$$
depending on the choice of chamber $\fC$ and polarization $\varepsilon$. For a component $Z\subset X^{A}$ and any $\gamma \in H^{\bullet}_{T}(Z)$ the stable envelope $\Gamma=\textrm{Stab}_{\fC,\varepsilon}(\gamma)$ is defined uniquely by the following properties:
\begin{itemize}
\item
$\textrm{supp}(\Gamma)\subset \textrm{Slope}_{\fC}(Z)$
\item $\left.\Gamma\right|_{Z}=\pm e(N_{-}) \cup \gamma$ with the sign chosen according to the polarization~$\varepsilon$.
\item $\deg_{A} \left.\Gamma\right|_{Z^{\prime}}< \textrm{codim}(Z^{\prime})/2$,  for any $Z^{\prime}> Z$
\end{itemize}
The proof of this theorem can be found in \cite{OM}. In addition, the description of $\textrm{Stab}_{\fC, \varepsilon}$
as the Lagrangian correspondence can be found there.

As we mentioned above, the choice of polarization is a formality corresponding to choice of signs. We will use symbol $\textrm{Stab}_{\fC}$ for the stable map meaning that some polarization $\varepsilon$ is chosen.

As an example consider the case $X=T^{\ast} {\mathbb{P}}^{n}$, with $A=({\mathbb{C}}^{\ast})^{n+1}$. The action of $A$ on $X$ is induced from canonical action on ${\mathbb{C}}^{n+1}$. Let $T= A \times {\mathbb{C}}^{\ast}$ where the additional factor ${\mathbb{C}}^{\ast}$ acts on the fibers of $T^{\ast}{\mathbb{P}}^{n}$ with a character $\hbar$.
The fixed set $X^{A}=\{p_0,...,p_n\}$, where $p_i$ are the points on ${\mathbb{P}}^{n}$ corresponding to the coordinate lines in ${\mathbb{C}}^{n+1}$.
The chambers in $a_{\mathbb{R}}$ are the standard chambers of $\frak{gl}(n+1)$. To describe the stable envelopes of the fixed  points $p_i\in X^A$, we choose $\fC$ to be the fundamental chamber of $\frak{gl}(n+1)$, such that is induces the natural ordering $p_n\succeq ... \succeq p_{0}$.

Next, note that the space $T^{\ast}{\mathbb{P}}^{n}$ equivariantly retracts to ${\mathbb{P}}^{n}$, and thus we have (for example see chapter 27.1 in \cite{PandLect}):
\be
\label{Pring}
H^{\bullet}_{T}(T^{\ast}{\mathbb{P}}^{n})={\mathbb{C}}[c,u_0,...,u_n,\hbar]/(c-u_0)(c-u_1)...(c-u_n)
\ee
where $c$ is the first Chern class of the tautological bundle ${\cal{O}}(-1)$ on ${\mathbb{P}}^{n}$ and $u_i$ are the equivariant parameters corresponding to the characters of $A$. Therefore, $\Stab_{\fC}(p_i)$ is a polynomial from  (\ref{Pring}) satisfying the conditions of the theorem above. We claim that:
\be
\label{sp}
\Stab_{\fC}(p_k)=\prod\limits_{i<k}\,(u_i-c-\hbar) \prod\limits_{i>k}\, (u_i-c)
\ee
One can check that this formula indeed enjoys all the necessary properties. The first condition $\textrm{supp}(\Gamma)\subset \textrm{Slope}_{\fC}(Z)$ means that
$$\left.\Stab_{\fC}(p_k)\right|_{p_i}:=\left.\Stab_{\fC}(p_k)\right|_{c=u_i}=0, \ \ \textrm{for}  \ \ k<i$$
what is obviously the case.
The third condition imply that:
$$
\deg_{u_1,...,u_n} \left.\Stab_{\fC}(p_k)\right|_{c=u_m}<n, \ \ \textrm{for} \ \ k>i
$$
what is also the case. Finally, the second condition means that the equivariant Euler class of $N_-$-component of the normal bundle to the point has the form:
$$
\left.e(N_-)\right|_{p_k}=\left.\Stab_{\fC}(p_k)\right|_{c=u_k}=\prod\limits_{i<k}\,(u_i-u_k-\hbar) \prod\limits_{i>k}\, (u_i-u_k)
$$
what is also true up to a sign which is just a choice of the polarization. Thus, by uniqueness of the stable envelope it proves (\ref{sp}).

\section{R-matrix \label{Rmsec}}

Assume that $X$ is as above, and the chamber $\fC$ together with polarization $\varepsilon$ are chosen. Then we have defined maps:
\be
\textrm{Stab}_{\fC,\varepsilon} : H_{T}^{\bullet}(X^{A})\rightarrow H_{T}^{\bullet}(X)
\ee
For a pair of two chambers, we can define the $R$-matrix as the following operator:
\be
\begin{array}{|c|}
\hline\\
R_{\fC^{\prime},\fC}=\textrm{Stab}^{-1}_{\fC^{\prime},\varepsilon} \circ \textrm{Stab}_{\fC,\varepsilon} \in \textrm{End}(H_{T}^{\bullet}(X^{A}))\otimes\mathbb{C}(\frak{t})\\
\\
\hline
\end{array}
\ee
where $\mathbb{C}(\frak{t})$ denotes rational functions on the Lie algebra $\frak{t}$ of  $T$.

\subsection{$R$-matrices for the flag varieties \label{flagR}}
An important case of the varieties satisfying conditions above are the Nakajima quiver varieties \cite{NakQuiv}, in particular the quiver varieties of  $A_{n}$ type:
$$
{\cal{A}}_n(w)=\coprod\limits_{v} {\cal{A}}_{n}(v,w)
$$
Here ${\cal{A}}_n (v,w)$ is the $A_{n}$ quiver varieties defined by the vectors of dimensions $v$ and $w$. We adopt the notations:
\be
w=\sum\limits_{k=1}^{n} w_k \delta_{k}, \ \ w=\sum\limits_{k=1}^{n} v_k \delta_{k}
\ee
for the dimension vectors with $\dim W_{i}=w_i$,  $\dim V_{i}=v_i$.
Note that fixing the dimensions $w$ impose certain condition for $v$. For example, as discussed in \cite{NakAle},\cite{NakQuiv}, if $w=m\delta_1$, then we have the following condition:
\be
w_1=m\geq v_1\geq v_2\geq...\geq v_{n}
\ee
and for the corresponding Nakajima variety we have:
\be
\label{flags}
{\cal{A}}_{n}(w)=\coprod\limits_{m\geq v_1\geq...\geq v_{n}}\, T^{\ast} F(v_1,v_2,...v_{n})
\ee
where $T^{\ast} F(v_1,v_2,...v_{n})$ is the cotangent bundle of the flag variety:
\be
F(v_1,v_2,...v_{n})=\{{\mathbb{C}}^m \supset V_{1} \supset V_{2} \supset...\supset V_{n}:  \dim V_{k}=v_k\}
\ee
Let us consider the action of the group
$$
G_{w}=\textrm{GL}(w_1)\times...\times \textrm{GL}(w_{n})
$$
on the representation of the quiver $W_1\oplus...\oplus W_{n}$. Let $A$ be a maximal torus of $G_w$. Let us consider the induced action of $A$  on ${\cal{A}}_n(w)$. This action clearly preserves the symplectic form.

We define $T=A\times {\mathbb{C}}^{\ast}$. The additional, one-dimension factor $ {\mathbb{C}}^{\ast}$ acts on ${\cal{A}}_n(w)$ as follows.
For $A_{n}$-quiver  ${\cal{A}}_n(v,w)$ consider the maps:
$$
X_k \in \textrm{Hom}(V_k,V_{k+1}), Y_k \in \textrm{Hom}(V_{k+1},V_k),   i_k \in \textrm{Hom}(W_{k},V_k),  j_k \in \textrm{Hom}(V_{k},W_k)
$$
Then the elements $z\in T/A\simeq {\mathbb{C}}^{\ast}$ act on the representation of the quiver as:
$$
z(X_k)=z X_k, \ \  z(Y_k)=Y_k, \ \  z(i_k)=z i_k, \ \  z(j_k)=j_k
$$
This induces an action of ${\mathbb{C}}^{\ast}$ that scales the symplectic form on ${\cal{A}}_n(w)$.
%Let $G\in G_{w}$ is the subgroup of elements $g=(g_1,g_2,...,g_{N-1})\in G_{w} $ where all components have the same determinant: $\det g_i=\det g_j$. The action of $G$ on the representation, induces an action on the total space ${\cal{N}}(w)$. Let $T$ is the maximal torus of $G$ and $A\subset T$ is the maximal torus of the subgroup:
%\be
%\textrm{SL}(w_1)\times...\times \textrm{SL}(w_{N-1})\subset \, \textrm{GL}(w_1)\times...\times \textrm{GL}_{w_{N-1}}
%\ee
%The action of $A\subset T$ on ${\cal{N}}(w)$ satisfies all necessary conditions: $A$-preserves the symplectic form and $T$ scales it. In particular, in the case
%$w=n \delta_1$, we have: $T=({\mathbb{C}}^{\ast})^N$ with the action on flags (\ref{flags}) induced from the standard scaling action on the $N$-dimensional %space, as in (\ref{scal}).

Let us consider the simplest Nakajima variety associated with the $A_{1}$-quiver for $w = 2 \delta_1$:
\be
\label{P1}
{\cal{A}}_1(2 \delta_1 )=\coprod\limits_{v\leq 2 } T^{\ast} F(v)
\ee
in this case $F(v)$ is the space of $v$-dimensional subspaces in ${\mathbb{C}}^2$. Thus, we have:
\be
{\cal{A}}_1(2 \delta_1 )=\textrm{pt}\coprod T^{\ast}{\mathbb{P}} \coprod \textrm{pt}
\ee
We first consider the stable map in cohomologies of $T^{\ast}{\mathbb{P}}$. As in section \ref{se}
the stable map has the form:
$$
\textrm{Stab}_{\fC_1}(p_1)=(u_2-c), \ \ \ \textrm{Stab}_{\fC_1}(p_2)=(u_1-c-\hbar)
$$
analogously for the opposite chamber we obtain:
$$
\textrm{Stab}_{\fC_2}(p_1)=(u_2-c-\hbar) \ \ \ \textrm{Stab}_{\fC_2}(p_2)=(u_1-c)
$$
Assume $i^{\ast}$ is the pullback on the cohomology induced by the inclusion of fixed points to $T^{\ast}{\mathbb{P}}$. In the equivariant situation the map $i^{\ast}$ is an isomorphism, so we can write:
$$
R_{\fC^{\prime},\fC}=\textrm{Stab}^{-1}_{\fC_2} \, \textrm{Stab}_{\fC_1}=\Big(i^{\ast} \circ \textrm{Stab}_{\fC_2} \Big)^{-1}  \Big(i^{\ast} \circ \textrm{Stab}_{\fC_1} \Big)
$$
The restriction of the class to a point $p_i$ amounts in the substitution $c=u_i$, such that in the ordered basis $p_1\preceq p_2$, the matrices of these operators take the form:
$$
i^{\ast}\circ \textrm{Stab}_{\fC_1}=\left(\begin{array}{cc} -u & -\hbar \\ 0 & u-\hbar \end{array}\right),  \ \ \ i^{\ast}\circ \textrm{Stab}_{\fC_2}=\left(\begin{array}{cc} -u-\hbar & 0 \\ -\hbar& u \end{array}\right)
$$
where $u=u_1-u_2$.
Therefore, $ T^{\ast}{\mathbb{P}}$ -part of $R$-matrix takes the form:
\be
\label{RmP}
R_{\fC_2,\fC_1}=\left(\begin{array}{cc}
\frac{u}{u+\hbar}& \frac{\hbar}{u+\hbar}\\
\frac{\hbar}{u+\hbar} & \frac{u}{u+\hbar}
\end{array}\right)
\ee
The whole $R$-matrix for (\ref{P1}):
\be
\label{Rmatf}
R=\left[ \begin {array}{cccc} 1&0&0&0\\ \noalign{\medskip}0&{\frac {u}{u+\hbar}}&{\frac {h}{u+\hbar}}&0\\ \noalign{\medskip}0&{\frac {h}{u+\hbar}}&{
\frac {u}{u+\hbar}}&0\\ \noalign{\medskip}0&0&0&1\end {array} \right]=\dfrac{u}{u+\hbar}\, \textrm{Id} +\dfrac{\hbar}{u+\hbar} \, \textrm{P}
\ee
where $\textrm{Id}$ and $\textrm{P}$ are the identity and permutation operators in $\mathbb{C}^2\otimes \mathbb{C}^2$. Note that this is the  well known  standard, rational  $\frak{gl}(2)$ $R$-matrix \cite{QISMfirst},\cite{Sklyanin}.

It is instructive to consider the same procedure for $A_{N-1}$-quiver with $w = 2 \delta_1$. In this case the variety has the form:
\be
\label{eeend}
{\cal{A}}_{N-1}(w)=\coprod\limits_{2\geq v_1\geq...\geq v_{N-1}}\, T^{\ast} F(v_1,v_2,...v_{N-1})
\ee
with the following dimension of the cohomology ring:
\be
\dim H^{\bullet} \Big(\coprod\limits_{2\geq v_1\geq...\geq v_{N-1}}\, T^{\ast} F(v_1,v_2,...v_{N-1}) \Big)= N^2
\ee
Calculation similar to one considered above leads to $R$-matrix of exactly the same form (\ref{Rmatf}):
\be
\label{Rsq}
\begin{array}{|c|}
\hline\\
 \ \ R=\dfrac{u}{u+\hbar} \, \textrm{Id} +\dfrac{\hbar}{u+\hbar}\,\textrm{P} \ \ \\
\\
\hline
\end{array}
\ee
where $\textrm{Id}$ and $\textrm{P}$ are the identity and permutation operators acting in ${\mathbb{C}}^{N}\otimes{\mathbb{C}}^{N}$. Note, that this $R$-matrix is the standard rational $\frak{gl}(N)$ $R$-matrix acting in the tensor square of the fundamental representation $V={\mathbb{C}}^{N}$. We see, that the cohomology of (\ref{eeend}) are identified with $V^{\otimes 2}$. Analogously, the cohomology of $A_{N-1}$ quiver variety with $w = m \delta_1$ can be identified with $V^{\otimes m}$. As we discuss in the next section, the  $R$-matrix for this case factorizes to certain products of the fundamental $R$-matrices (\ref{Rsq}).

In the general case of $A_{N-1}$ with $w = \sum m_{i} \delta_i$ we have the identification:
\be
H^{\bullet} \Big({\cal{A}}_{N-1}(w )\Big) =\bigotimes\limits_{k=1}^{N-1} \left({\bigwedge}^{k} V\right)^{\otimes m_k}
\ee
and the resulting $R$-matrix factorizes to a product of $R$-matrices for the other fundamental representations  $\bigwedge^{m_i} V \otimes \bigwedge^{m_j} V$. The last ones, can be computed from (\ref{Rsq}) by means of the fusion procedure described in section~\ref{FPR}.  For a general quivers associated to each Dynkin diagram the resulting $R$-matrices are the rational $R$-matrices for corresponding Lie algebras. The formal proof of this fact was given in \cite{Michael} and is based on the results of \cite{Nak6} and  \cite{Varagnolo}.

\subsection{Factorization of $R$-matrices \label{FactR}}
Let us consider the case $A={\mathbb{C}}^{\ast}$, then its real Lie algebra defined by (\ref{rlie}) is one-dimensional $a_{{\mathbb{R}}}\simeq {\mathbb{R}}$. Thus, $0 \in a_{{\mathbb{R}}}$ is the only codimension one subspace in it, and the chamber decomposition of $a_{{\mathbb{R}}}$ in any case is of the form:
$$
a_{{\mathbb{R}}}\setminus \{0\}=\{u>0\} \cup \{u<0\}
$$
As we have only two chambers $\fC=\{u>0\}$, $\fC^{\,\prime}=\{u<0\}$ there is only one nontrivial  $R$-matrix we can construct:
$$
R(u)=\textrm{Stab}^{-1}_{\fC^{\,\prime}} \circ \textrm{Stab}_{\fC}
$$
Let us consider what can happen when the torus $A$ gets larger, i.e. we have  an inclusion $A \hookrightarrow \hat{A}$ to a bigger torus $\hat{A}$ acting on
$X$. If $\dim \hat{A} \geq 2$, we can have many, and possibly even infinite number of chambers $\fC_{\,i}$ in the corresponding Lie algebra $\hat{a}_{{\mathbb{R}}}$.
The initial Lie algebra $a_{{\mathbb{R}}}\simeq{\mathbb{R}}$ is some line in $\hat{a}_{{\mathbb{R}}}$. We can assume that $a_{{\mathbb{R}}}^{+}=\{u>0\}$ is
in one of these new chambers $\fC_{+}$ and $a_{{\mathbb{R}}}^{-}$ is in the opposite one $\fC_{-}=-\fC_{+}$. As in the new Lie algebra there are many chambers, to pass from the chamber $\fC_{+}$ to $\fC_{-}$ in $\hat{a}_{{\mathbb{R}}}$ we may have to cross some walls separating chambers $\fC_{+} ,\fC_{1}$; $\fC_{1},\fC_{2}$; ... ;$\fC_{m},\fC_{-}$. And taking into account that the maps $\textrm{Stab}_{\fC_{i}}$ are isomorphisms we can rewrite the initial $R$-matrix in the form:
$$
\begin{array}{l}
R(u)=\textrm{Stab}_{\fC_{\,+}}^{-1} \circ \textrm{Stab}_{\fC_{\,-}}=\\
\\
=\Big(\textrm{Stab}^{-1}_{\fC_{+}} \circ \textrm{Stab}_{\fC_{1}} \Big)\circ\Big( \textrm{Stab}_{\fC_{1}}^{-1} \circ \textrm{Stab}_{\fC_{2}}\Big)...\circ \Big(\textrm{Stab}_{\fC_{m}}^{-1}\circ \textrm{Stab}_{\fC_{-}}\Big)=
\\
\\
=R_{\fC_{+},\fC_{1}}  R_{\fC_{1},\fC_{2}}...R_{\fC_{m},\fC_{-}}.
\end{array}
$$
Thus, after enlarging the torus, the $R$-matrix factorizes to a product $R$-matrices corresponding to the new adjacent chambers. As any two adjacent chambers are separated by a wall corresponding to some root, we call them  \textit{root $R$-matrices}.
\begin{figure}
\begin{center}
%{\includegraphics*[natwidth=20cm,natheight=7cm,scale=0.7]{block1.eps}}
\unitlength 1mm % = 2.845pt
\linethickness{0.4pt}
\ifx\plotpoint\undefined\newsavebox{\plotpoint}\fi % GNUPLOT compatibility
\begin{picture}(100,30)(0,2)
\put(48,30){\vector(-1,0){25}}
\put(48,30){\vector(1,0){25}}
\put(48,30){\vector(2,3){13}}
\put(48,30){\vector(-2,3){13}}
\put(48,30){\vector(-2,-3){13}}
\put(48,30){\vector(2,-3){13}}
\put(75,29){$\textrm{ker}_{\alpha_{1}}$}
\put(30,52){$\textrm{ker}_{\alpha_{2}}$}
\put(60,52){$\textrm{ker}_{\alpha_{3}}$}
\put(48,30){\line(-2,1){20}}
\put(48,30){\line(2,-1){20}}
\put(70,18){$a_{{\mathbb{R}}}^{+}$}
\put(21,42){$a_{{\mathbb{R}}}^{-}$}
\end{picture}
\caption{\label{fg1} \footnotesize{ The chamber decomposition for $T^{\ast}{\mathbb{P}}^2$
}}\label{4block}
\end{center}\end{figure}
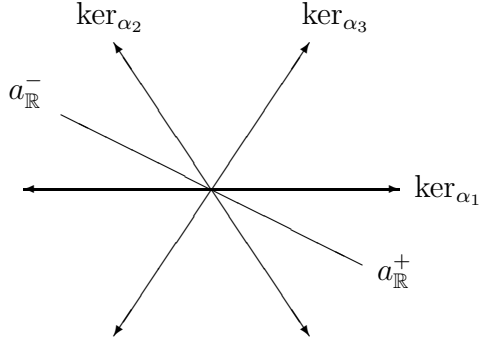

Let us return to the example considered in the section \ref{flagR}. For instance, consider the action of $\hat{A}=({\mathbb{C}}^{\ast})^2$ on $T^{\ast}{\mathbb{P}}^2$, defined in the homogenous coordinates~by:
\be
(z_1,z_2,z_3)\cdot ( x:y:z )=( z_1 x:z_2 y:z_3 z ), \ \ \textrm{ with}  \ \ z_1 z_2 z_3=1
\ee
The real Lie algebra $\hat{a}_{{\mathbb{R}}}$ is divided into six chambers by the walls corresponding to the roots $\alpha_1=u_1-u_2$, $\alpha_2=u_2-u_3$, $\alpha_3=u_1-u_3$, as in fig.~\ref{fg1}. We choose one-dimensional torus $A\subset \hat{A}$, such that its real Lie algebra ${a}_{{\mathbb{R}}}$ is as in fig.~\ref{fg1}.  To pass from ${a}_{{\mathbb{R}}}^{+}$ to ${a}_{{\mathbb{R}}}^{-}$ in $A\subset \hat{A}$ we have to cross three walls. Moreover, we have two ways to do it. Thus we obtain two different factorizations of the $R$-matrix for the torus $A$:
\be
\label{YB}
R_{12}(u_{12})R_{13}(u_{13})R_{23}(u_{23})=R_{23}(u_{23})R_{13}(u_{13})R_{12}(u_{12}), \ \ u_{i j}=u_i-u_j
\ee
Let us consider how each of the root matrices $R_{i j}(u_i-u_j)$ acts. In this particular case, the fixed sets of $A$ are three points on the base  ${\mathbb{P}}^2$ of the bundle  $T^{\ast}{\mathbb{P}}^2$ corresponding to the coordinate lines. As we discuss in section~\ref{cdec}, the walls correspond to the coweights for which the fixed sets get larger. For example, the wall for $\alpha_1=u_1-u_2$ corresponds to the coweights fixing the equivariant subspace
$T^{\ast}{\mathbb{P}}\subset T^{\ast}{\mathbb{P}}^2$ connecting the fixed points $p_1$ and $p_2$.  Now, if  we denote by $\fC_1$ and $\fC_2$ the chambers separated by the wall $u_1-u_2$. The action of $\textrm{Stab}_{\fC_1}$ and $\textrm{Stab}_{\fC_2}$ in the basis of the fixed points differs only on the components corresponding to the classes of the fixed points $p_1$ and $p_2$ such that:
\be
R_{1 2}(u_1-u_2)=R(u_1-u_2)\otimes 1
\ee
meaning that on the class of the third fixed point the action of this $R$-matrix is trivial $R_{1 2}(u_1-u_2)p_3=p_3$.

The $R$-matrix $R(u_1-u_2)$ acts in the basis of classes $p_1$ and $p_2$ as the $R$-matrix for $T^{\ast}{\mathbb{P}}^{1}$, explicitly as given by (\ref{RmP}). Similar consideration applies  to the other root matrices. Thus, the above formula (\ref{YB}) is exactly the Yang-Baxter relation for (\ref{RmP}).

In general, the root $R$-matrices satisfy the Yang-Baxter equation when, for example, the corresponding roots form the standard root system of~$\frak{gl}(N)$. In this case the walls are the standard walls separating the Weyl chambers. The intersection of any two walls is a facet surrounded by six Weyl chambers in the Lie algebra. To pass from one of them into the opposite we have to cross three walls in two possible ways, what give us factorization~(\ref{YB}). This is the case for any $A_{N-1}$ quiver variety, what proves, in particular, the Yang-Baxter relation for the $R$-matrix (\ref{Rsq}) considered above.

We can now continue example of the previous section: $A_{N-1}$-quiver with~$w = m \delta_1$. In this case we have:
\be
\label{eee}
{\cal{A}}_{N-1}(w)=\coprod\limits_{m \geq v_1\geq...\geq v_{N-1}}\, T^{\ast} F(v_1,v_2,...v_{N-1})
\ee
The dimension of its cohomology space is:
$$
\dim H_{T}^{\bullet}\Big({\cal{A}}_{N-1}(m \delta_1)\Big)=N^m
$$
and this space can be identified with the tensor power of fundamental representation for $\frak{gl}(N)$:
\be
\label{tpr}
H_{T}^{\bullet}\Big({\cal{A}}_{N-1}(m \delta_1)\Big)\simeq V^{\otimes m}
\ee
The roots for ${\cal{A}}_{N-1}(w)$ form the standard $\frak{gl}(N)$ root system.  The fixed set ${\cal{A}}_{N-1}^A(m \delta_1)$ is a disjoint union of
$N^m$ points (they are the flags formed by coordinate subspaces ). The walls between the Weyl chambers correspond to some fixed $T^{\ast}{\mathbb{P}}$ in ${\cal{A}}_{N-1}(w)$. Thus, the resulting $R$-matrix factorizes to a product of (\ref{Rsq}). Explicitly, the factorization procedure gives:
\be
\label{prdfor}
R=\prod\limits_{i=1}^{n-1} \prod\limits_{j=i+1}^{n} \, R_{ij}(u_i-u_j)
\ee
where as usual, $R_{ij}(u_i-u_j)$ stands for $R$-matrix (\ref{Rsq}) acting in $i$-th and $j$-th component of the tensor product (\ref{tpr}), and we imply that the product is ordered from the right to the left such that, for example, for $n=3$ we will have:
$$
R=R_{23}(u_2-u_3) R_{13}(u_1-u_3) R_{12}(u_1-u_2)
$$

In general, the root $R$-matrix for a root $\alpha$ is the $R$-matrix corresponding to the subspace fixed by cocharacters on the wall $\ker_{\alpha}$.
In some cases these subspaces have simpler structure, such that the corresponding root $R$-matrices are known. In the next section we show that the $R$-matrix for the instanton moduli space factorizes to a infinite product of $R$-matrices of $A_{\infty}$ type, which will be computed explicitly in section \ref{comp}.

\subsection{$R$-matrix for the instanton moduli space \label{instR}}
Now we return to the ADHM description of the instanton moduli space outlined in section \ref{embd}. As we noted there, the fixed set of torus $A$ has the form (\ref{emb}). The roots  form the standard $\frak{gl}(r)$ root system.  Thus, by factorization, described in the previous section, for the $R$-matrix we have:
\be
{\cal{R}}=\prod\limits_{i=1}^{r-1} \prod\limits_{j=i+1}^{r} \, {\cal{R}}_{ij}(u_i-u_j)
\ee
where ${\cal{R}}_{ij}(u_i-u_j)$ is the corresponding $R$-matrix acting in the $i$-th and $j$-th tensor component of:
\be
H^{\bullet}_{T}\Big(  \coprod\limits_{n=0}^{\infty} {\cal{M}}^{A}(n,r)\Big)\simeq{\cal{F}}^{\otimes r} \otimes\mathbb{C}[\frak{t}]
\ee
for example, ${\cal{R}}_{1,2}(u)={\cal{R}}(u)\otimes 1\otimes ...\otimes 1$. Therefore, as in the case of flags, the $R$-matrices for general $r$
can be written as products of $R$-matrices ${\cal{R}}(u)$ for $r=2$. In the following text we focus on the case $r=2$.

Let us consider the case $r=2$ in details. First, we choose the subgroup of the torus $A$ explicitly given by elements:
\be
\label{chm}
g_w=\left(\begin{array}{cc}
z & 0\\
0 & 1
\end{array}\right)
\ee
The fixed set is of the form:
\be
{\cal{M}}^{A}(2)=\coprod\limits_{n=0}^{\infty} {\cal{M}}^{A}(n,2)=\coprod\limits_{n=0}^{\infty} \Big(\coprod\limits_{n_1+n_2=n}\, \textrm{Hilb}_{n_1}\times \textrm{Hilb}_{n_2}\Big)
\ee
The chosen torus is one-dimensional and have two chambers $u\lessgtr 0$. The corresponding $R$-matrix acts in:
\be
H_{T}^{\bullet}\left({\cal{M}}^{A}(2)\right) \otimes\mathbb{C}[\frak{t}]={\cal{F}}\otimes {\cal{F}} \otimes\mathbb{C}[\frak{t}]
\ee
The main idea of computing this $R$-matrix is to consider the inclusion of this torus to a bigger two-dimensional one. Then, to pass from $u>0$ to $u<0$ in the new Lie algebra we will have to cross infinitely many walls. This gives the factorization of the $R$-matrix into the infinite product of certain operators which will be simple to compute explicitly. The additional dimension comes from the subtorus of $\textrm{GL}(2)$ in (\ref{gr}) given by elements of the form:
$$
\left(\begin{array}{cc}
w & 0\\
0 & w^{-1}
\end{array}\right)\in \textrm{GL}(2)
$$
%such that the we have the following tree dimensional torus acting acting on moduli space
%\be
%\left(\begin{array}{cc}
%z & 0\\
%0 & 1
%\end{array}\right)\times \left(\begin{array}{cc}
%w_1 & 0\\
%0 & w_2
%\end{array}\right) = T
%\ee
%We denote by $(u,t_1,t_2)$ the character of $T$ corresponding to the representation:
%\be
%\left(\begin{array}{cc}
%z^u & 0\\
%0 & 1
%\end{array}\right)\times \left(\begin{array}{cc}
%z^{t_1} & 0\\
%0 & z^{t_2}
%\end{array}\right) \subset T
%\ee
%The two dimensional subtorus $\hat{A}\subset T$ fixing the symplectic form on ${\cal{M}}(n,2)$ is defined by condition
%\be
%\label{hb}
%\hbar=t_1+t_2=0.
%\ee
The enlargement $\hat{A}$ of the torus of $A$ is given explicitly by the elements of the form:
\be
\label{subt}
\hat{A}=\left\{\left(\begin{array}{cc}
z & 0\\
0 & 1
\end{array}\right)\times \left(\begin{array}{cc}
w & 0\\
0 & w^{-1}
\end{array}\right) \right\} \subset T
\ee
The torus $\hat{A}$ fixes the symplectic form on ${\cal{M}}(n,2)$ and thus satisfies all necessary conditions. Let us consider the fixed set ${\cal{M}}^{\hat{A}}(n,2)$. Under the action of $\hat{A}$ the representation of the quiver (\ref{adhm}) splits to the irreducible ones:
$$
V=\bigoplus\limits_{k,m} \, V^{(k,m)}, \ \ \ W=\bigoplus\limits_{k,m} \, W^{(k,m)}, \ \ \left.g \right|_{V^{(k,m)}}=\left.g \right|_{W^{(k,m)}}=z^k w^m
$$
If $(B_{1},B_{2},i,j)$ represents a fixed point on (\ref{adhm}), then we have:
\be \label{stp}
j=g_w j g_v^{-1}, \ \ i= g_v i g_w^{-1}, \ \  g_v B_1 g_v^{-1}= w B_1, \ \  g_v B_2 g_v^{-1}= w^{-1} B_2
\ee
The first two equations mean that the corresponding operators split, such that the only nontrivial components of $i$ and $j$ are:
$$
i=i_0\oplus i_1, \ \ j=j_0\oplus j_1, \ \ i_k\in \textrm{Hom}(W^{(k,0)},V^{(k,0)}), \ \ j_k\in \textrm{Hom}(V^{(k,0)},W^{(k,0)})
$$
where $k=0,1$,
such that $W=W^{(0,0)}\oplus W^{(1,0)}$ with $\dim W^{(0,0)}=\dim W^{(1,0)}=1$.
The second pair of equations in (\ref{stp}) means that the matrices
$B_i$ have the following block form:
$$
B_{i}=\bigoplus\limits_{m\in {\mathbb{Z}},\,k=0,1} B_{i}^{(k,m)}, \ \ B_{1}^{(k,m)}\in \textrm{Hom}(V^{(k,m)},V^{(k,m+1)}), \ \ B_{2}^{(k,m)}\in \textrm{Hom}(V^{k,m-1},V^{k,m})
$$
%We note here, that the matrices $B_i$ splits completely, i.e. do not have components in $\textrm{Hom}(V^{k_1,m-1},V^{k_2,m})$ with $k_1\neq k_2$.
The data $(V^{(k,m)}, W^{(k,m)}, i_k, j_k)$  for each value $k=0, 1$ defines a representation of $A_{\infty}$ quiver  with $\dim W_{m}=\delta_{m,0}$. Therefore, the fixed set is of the form:
\be
\label{fixs}
{\cal{M}}^{\hat{A}}(n,2)= \coprod\limits_{\sum_i p_i +\sum_i q_i=n }{\cal{A}}_{\infty}(p_i, \delta_0) \times {\cal{A}}_{\infty}(q_i, \delta_0)
\ee
where ${\cal{A}}_{\infty}(p_i, \delta_0)$ is the $A_{\infty}$ quiver variety corresponding to the vectors of dimensions $\dim V_i=p_i$, $\dim W_i =\delta_{i,0}$.
We also note  here that the set $\coprod\limits_{\sum_i p_i=m}{\cal{A}}_{\infty}(p_i, \delta_0)$ is the disjoint union of points $p_\lambda$ labeled by the partitions with $|\lambda|=m$. Thus, for the fixed set (\ref{fixs}) we have:
$${\cal{M}}^{\hat{A}}(n,2)=\coprod\limits_{|\lambda|+|\mu|=n}\, p_\lambda \times p_\mu$$
%This is, of course, the standard fixed points on the Hilbert schemes:
%$$
%{\cal{M}}^{\hat{A}}(n,2)=\Big({\cal{M}}^{A}(n,2) \Big)^{\hat{A}}=\Big(\coprod\limits_{n_1+n_2=n} \textrm{Hilb}_{n_1}\times \textrm{Hilb}_{n_2}  %\Big)^{\hat{A}}
%$$
Next, we need to identify the chambers in the Lie algebra of $\hat{A}$.
 %Let denote by $u$ the character of torus $A \in \hat{A}$ and by $t_1$ the characters of %$\hat{A}/A$ such that $u, t_1$ generate $c(\hat{A})$. We are going to show that the the Lie algebra $\hat{a}_{\mathbb{R}}$ is separated into infinitely many %chambers by walls defines as lines $u+m t_1=0$ for $m\in \mathbb{Z}$.
Let us fix a cocharacter, defining some one-dimensional subtorus $Q\in\hat{A}$ consisting of the elements:
\be
Q=\left\{\left(\begin{array}{cc}
z^{u} & 0\\
0 & 1
\end{array}\right)\times \left(\begin{array}{cc}
z^{t_1} & 0\\
0 & z^{-t_1}
\end{array}\right)\right\}
\ee
We need to understand for which values of $u$ and $t_1$  the fixed set ${\cal{M}}^{Q}(n,2)$ gets larger than ${\cal{M}}^{\hat{A}}(n,2)$. These values correspond to the cocharacters lying on the walls.

%%%%%%%%%%%%%%%%%%%%%%%%%%%%%%%%%%%%%%%%%%%%%%%%%%%%%%%%%%%%%%%%%%%%%%%%%%%%%%%%%%
The action of this torus  on $n\times n$ matrices $B_i$ is of the form: $B_1\rightarrow z^{t_1} B_1, B_2\rightarrow z^{-t_1} B_2$ . This point is fixed, in the quotient (\ref{adhm}) if there exists $g_v\in \textrm{GL}(n)$ such that:
\be
z^{t_1} B_1=g_v B_1 g_v^{-1}, \ \ z^{-t_1} B_2=g_v B_2 g_v^{-1}
\ee
If $V=\bigoplus\limits_{k\in {\mathbb{Z}}} V^{(k)}$ is the weight decomposition of $V$ under the action of $Q$, such that $g_v$ on each component $V^{(k)}$ acts diagonally with weights $k$, i.e. by multiplication on $z^{k\, t_1}$ then the last equations imply that the matrices $B_{i}$ have a block form:
\be
\label{Bs}
B_{i}=\bigoplus\limits_{k\in {\mathbb{Z}}} B_{i}^{(k)},  \ \ B_{1}^{(k)}\in \textrm{Hom}(V^{(k)},V^{(k+1)}), \ \ B_{2}^{(k)}\in \textrm{Hom}(V^{(k-1)},V^{(k)})
\ee
and all other components of $B_i$ vanish.
Similarly $j$ is $n\times 2$ matrix which represents the fixed point if the following condition holds:
\be
g_{w} j g_{v}^{-1}=j
\ee
it is possible only if $j$ splits:
\be
\label{js}
j=j_0\oplus j_m, \ \ \ j_0 \in \textrm{Hom}(V^{(0)}, W^{(0)}), \ \ \ j_m \in \textrm{Hom}(V^{(m)}, W^{(u)})
\ee
and $u=- m t_{1}$ for some $m \in {\mathbb{Z}}$. The same is true for $i$:
\be
\label{is}
i=i_0\oplus i_m, \ \ \ i_0 \in \textrm{Hom}(W^{(0)}, V^{(0)}), \ \ \ i_m \in \textrm{Hom}(W^{(u)}, V^{(m)})
\ee
where $W=W^{0}\oplus W^{u}$, such that $\left.g_w\right|_{ W^{0}}=1$ ,$\left.g_w\right|_{W^{u}}=z^u$, and $\dim W^{0}=\dim W^{u}=1$.
Now, note that  (\ref{Bs}), (\ref{js}) and (\ref{is}) is nothing but data defining $A_{\infty}$ quiver, with representations $V^{(k)}$ and $W^{(k)}$.
The auxiliary spaces $W^{(k)}$ have dimensions $W^{(k)}=1$ for $k=0$ and $k=m$ and $0$ otherwise. Thus, we obtain:
\be
\label{fs}
{\cal{M}}^{Q}(n,2)=\coprod\limits_{\sum\,{n_i}=n}{\cal{A}}_{\infty}(n_i,\delta_0+\delta_m)
\ee
Now, we can use the factorization procedure, to describe the instanton $R$-matrix.
In the two-dimensional  torus $\hat{A}$ with characters $u, t_1$, to pass from the original chamber $u>0$ to $u<0$ in $A$ we need to cross infinitely many walls defined by the lines $u+ m t_1=0$ for $m\in {\mathbb{Z}}$. The factorization procedure, in this case gives:
\be
\label{prodFor}
\begin{array}{|c|}
\hline
\\
 \ \ \ {\cal{R}}(u)=\prod\limits_{m\in {\mathbb{Z}}}^{\rightarrow} \, R_{m}(u+m t_1) \ \ \  \\
\\
\hline
\end{array}
\ee
where the  root $R$-matrices $R_{m}(u)$ are the $R$-matrices for the fixed sets (\ref{fs}):
$$
R_{m}(u) \in \textrm{End}\left(  H^{\bullet}_{T}\Big(\coprod\limits_{n=0}^{\infty}\coprod\limits_{\sum\,{n_i}=n}{\cal{A}}_{\infty}(n_i,\delta_0+\delta_m)\Big) \right)
$$
This $R$-matrices can be identified with the so called universal $R$-matrices for the half-infinite wedge product of the fundamental representation of $\frak{gl}_\infty$. They can be computed explicitly. In fact, it is the generalization of $R$-matrix (\ref{Rsq}) that we have already computed for flag manifolds.  We describe these calculations in the next section.
\section{Computation of $R$-matrices \label{comp}}
The main result of the previous section is the product formula for instanton $R$-matrix (\ref{prodFor}). In this section we derive explicit formula for its factors $R_{m}(u)$ which are the $R$-matrices associated with $A_{\infty}$ quiver variety (\ref{fs}).
\subsection{$A_{n}$ $R$-matrices}
We start by a short outline of the fundamental representations of $\frak{gl}(N)$. Let  $\omega_{i},\ \ i=1...N-1$ be the set of fundamental weights of $\frak{gl}(N)$. Denote by $L(\omega_{i})$ the  $i$-th fundamental representation of  $\frak{gl}(N)$. By  Weyl dimension formula we obtain:
$$
\dim L(\omega_{i})=\dfrac{N!}{i! (N-i)!}
$$
The first $N$-dimensional module $L(\omega_{1})$ is the standard representation of $\frak{gl}(N)$ in $V={\mathbb{C}}^N$. If $e_{i}, \ \ i=1...N$ is the basis in $V$ and $E_{i j}$ are the standard $\frak{gl}(N)$ basis then this representation is given explicitly by:
$$
E_{i j} e_{k} = \delta_{jk} e_{i}
$$
Formally, this action is given by the following differential operators:
\be
\label{do}
E_{i j}=e_{i} \frac{\partial}{\partial e_{j}}
\ee
The other fundamental representations $L(\omega_{k})$ are identified with the wedge powers of the first one:
\be
\nonumber
L(\omega_{k})=\bigwedge^{k} V, \ \ \textrm{basis}\,L(\omega_{k}) =\{\, e_{i_1}\wedge e_{i_{2}}\wedge...\wedge e_{i_k}, \, 1\leq i_{1}<i_{2}<..<i_{k}\leq N\, \}
\ee
The action of the generators $E_{ij}$ on the basis elements of $L(\omega_{k})$ is given by the same differential operators (\ref{do}), where by multiplication on $e_{i}$ we understand the wedge product $e_{i}\,\wedge$ and by derivative the corresponding adjoint operator. Note, that these operators $e_{i}$ and $\frac{\partial}{\partial e_{i}}$ are anticommuting. It motivates to introduce the following \textit{fermions}:
\be
\label{ac}
\psi_{k}=e_{k} \wedge, \ \ \ \psi_{k}^{\ast}=\frac{\partial}{\partial e_{i}}
\ee
which satisfy the standard Clifford algebra relations:
\be
\label{cl}
\psi_{i} \psi_{j} +\psi_{j} \psi_{i}=0, \ \ \psi_{i}^{\ast} \psi_{j}^{\ast} +\psi_{j}^{\ast} \psi_{i}^{\ast}=0, \ \ \psi_{i}^{\ast}  \psi_{j}+ \psi_{j} \psi_{i}^{\ast} =\delta_{i j}
\ee
In this way, all fundamental representations of $\frak{gl}(N)$ can be treated at one go - the sum of them is the standard representation of Clifford algebra (\ref{cl}):
\be
\Lambda^{\bullet} V=\bigoplus\limits_{k=0}^{N} L(\omega_{k})
\ee
(cases $k=0,N$ correspond to trivial representations) and the action of $\frak{gl}(N)$ on each component is given by:
\be
\label{Ephi}
E_{ij}=\psi_{i} \psi_{j}^{\ast}
\ee
Let
$${\cal{A}}_{n}(w)=\coprod\limits_{v}{\cal{A}}_{n}(v,w) $$
be the $A_{n}$ quiver variety defined by the dimension vector $w$. The cohomology of this space carries a natural structure of $\frak{gl}(n+1)$ modules \cite{NakAle}, \cite{NakQuiv}. In particular, the fundamental representations $L(\omega_{k})$ correspond to the choice of the dimension vector $w=\delta_k$ i.e. all but one dimensions vanish: $\dim w_{i}=\delta_{i,k}$. In this case we have
$$
L(\omega_{k}) \simeq H^{\bullet} \Big(\coprod\limits_{v}{\cal{A}}_{n}(v,\delta_k)  \Big)
$$
For the general $A_{n}$-quiver variety, with $w=\sum_{k} m_k \delta_k$, i.e. for dimensions $\dim w_k =m_k$  we have:
\be
\label{gq}
H^{\bullet} \Big({\cal{A}}_{n}(w) \Big)\simeq L^{\otimes m_{1}}(\omega_{1})\otimes...\otimes L^{\otimes m_{n}}(\omega_{n})
\ee

Let $R(k_1,k_2)$ be the $R$-matrix corresponding to the quiver variety with dimension $w=\delta_{k_1}+\delta_{k_2}$. This $R$-matrix acts in the tensor product of two fundamental representations:
\be
\label{ccc}
H^{\bullet} \Big(\coprod\limits_{v}\, {\cal{A}}_n(\delta_{k_1}+\delta_{k_2}) \Big)=L(\omega_{k_1})\otimes L(\omega_{k_2})= {\bigwedge}^{k_1} V \otimes {\bigwedge}^{k_2} V
\ee
The $R$-matrix for the general $A_n$ quiver variety (\ref{gq}) are  always given by certain products of "elementary" $R$-matrices $R(k_1,k_2)$ as it was shown in the section~\ref{FactR}  for the case $w=m\delta_1$. Explicitly, this product is of same form~(\ref{prdfor}):
$$
R=\prod\limits_{i=1}^{m_1+...+m_{N-1}} \prod\limits_{j=i+1}^{m_1+...+m_{N-1}} \, R_{ij}(u_i-u_j)
$$
where $R_{i,j}$ are the $R(k_1,k_2)$-matrices acting in the $i$-th and $j$-th tensor component of (\ref{gq}).
Therefore, to describe all $R$-matrices for $A_n$ quivers we need to know $R(k_1,k_2)$. This elementary building blocks of the theory can be computed by the so called \textit{fusion procedure}.
\subsection{Fusion procedure \label{FPR}}

The fusion procedure allows constructing the $R$-matrices for general representations of $\frak{gl}(N)$ from the elementary building block - the $R$-matrix for the first fundamental representation. The  $R$-matrix for the first fundamental representation is the element $\textrm{End}({\mathbb{C}}^N\otimes {\mathbb{C}}^N)$ that was computed in  (\ref{Rsq}):
\be
R(u)=\dfrac{u}{u+\hbar}\,\textrm{ Id} +\dfrac{\hbar}{u+\hbar}\,\textrm{ P}
\ee
where $\textrm{Id}$ is the identity operator and $\textrm{P}$ is the permutation: $\textrm{P}(e_{i} \otimes e_{j}) =e_{j} \otimes e_{i}$.
This operator satisfies the quantum Yang-Baxter equation in $\mathbb{C}^N\otimes\mathbb{C}^N\otimes\mathbb{C}^N$:
\be
R_{12}(u) R_{13}( u+ v) R_{23}(v)=R_{23}(v) R_{13}( u+v) R_{12}(u)
\ee
where, as usual, $R_{n m}(u)$ is the operator acting in the $n$-th and $m$-th spaces and is identity in the third one. The fusion procedure allows producing more complicated $R$-matrices for a tenor product of two general representations of  $\frak{gl}(N)$: $R^{\mu,\nu}(u) \in \textrm{End} (V^{\mu}\otimes V^{\nu})$, satisfying the quantum Yang-Baxter equation in $V^{\lambda}\otimes V^{\mu}\otimes V^{\nu} $:
\be
R^{\lambda, \mu }_{12}(u) R^{\lambda, \nu}_{13}( u+ v) R^{\mu,\nu}_{23}(v)=R^{\mu,\nu}_{23}(v) R^{\lambda, \nu}_{13}( u+ v) R^{\lambda, \mu }_{12}(u).
\ee
We outline this procedure here. Let $V^{\lambda}$ be a $\frak{gl}(N)$ representation corresponding to the Young diagram $\lambda=(\lambda_1,\lambda_2,...,\lambda_{m})$ with  $n=\sum_{i} \lambda_{i}$ boxes. Let $P_{\lambda}: V^{\otimes n} \rightarrow V^{\lambda}$ be the Young projector from the $n$-th tensor degree of the fundamental representation to  $V^{\lambda}$. To the box with coordinates $(i,j)$ in the Young diagram we attach the following number:
\be
s_{(i j)}= u + (i-j) \hbar
\ee
%Such that for example in the case $\lambda=[3,2,2,1]$ we have Fig.???.
Enumerating all boxes in alphabetical order (i.e. from the left to the right in the first row, then from the left to the right in the second one, and so on), one gets the ordered sequence of numbers $(s_{1},s_{2},...,s_{n})$.

By construction, (see \cite{Zabrodin} for details) the $R$-matrix acting in $V^{\lambda}\otimes V$ is given by:
\be
\label{fp}
R^{\lambda}(u)=(P_{\lambda}\otimes 1 ) \, R_{n 0}(s_n ) \otimes ... \otimes R_{2 0}(s_2 ) \otimes R_{1 0}(s_1)\, (P_{\lambda}\otimes 1 )
\ee
The $R$-matrices $R_{n 0}(s_n )$ are operators acting in $V^{\otimes (n+1)}$ with one common "auxiliary"  space $V_{0}=\mathbb{C}^N$. The equation (\ref{fp}) implies that the operator $A_{\lambda}=R_{n 0}(s_n ) \otimes ... \otimes R_{1 0}(s_1) \in \textrm{End} (V^{\otimes (n+1)})$ preserves the subspace
$V^{\lambda}\otimes V \in V^{\otimes (n+1)}$, i.e. if $v\in V^{\lambda}\otimes V $ then $A_{\lambda}(v)\in V^{\lambda}\otimes V$ and the restriction of the operator $A_{\lambda}$ to the subspace $V^{\lambda}\otimes V$  gives the $R$-matrix for $ V^{\lambda}\otimes V$.

Applying the same procedure for $R^{\lambda}(u)\in  \textrm{End} (V^{\lambda}\otimes V)$ and taking the first space $V^{\lambda}$ as "auxiliary" we can construct the $R$-matrix $R^{\lambda,\nu}(u)$ for $V^{\lambda}\otimes~V^{\nu}$.

For  skewsymmetric representations $\bigwedge^{n} V$, when the corresponding partition is given by $\lambda=[\underbrace{1,1,...,1]}_{n}$
the fusion procedure gives the following result: the $R(n_1,n_2)$-matrix for $\bigwedge^{n_{1}} V \otimes \bigwedge^{n_{2}} V$ is given by the operator:
\be
\label{Rnm}
R(u)={\prod_{m=0}^{n_1-1}}^{\leftarrow}{\prod_{k=1}^{n_2}}^{\leftarrow} R_{n_1-m, n_1+k} \Big(u + (m+k-1) \hbar \Big)
\ee
The arrows mean that the products are ordered from the right to the left (we omit the projectors from (\ref{fp}) for simplicitly).  For example in the case $\bigwedge^{2} V \otimes \bigwedge^{3} V$, when $n_1=2$ and $n_2=3$  we have the following product:
\be
\label{R23}
R=R_{1 5}(u+3\hbar) R_{1 4}(u +2\hbar)  R_{1 3}(u +\hbar) R_{2 5}(u +2\hbar) R_{2 4}(u +\hbar ) R_{2 3}(u)
\ee
The operator (\ref{Rnm}) has invariant subspace $\bigwedge^{n_1} V \otimes \bigwedge ^{n_2} V \subset V^{\otimes n_1 }\otimes V^{\otimes n_2 } $. The restriction of the operator to this subspace gives the operator we need. The resulting $R$-matrix satisfies the Yang-Baxter equation what easily follows from the fact that each multiple in (\ref{Rnm}) satisfies YBE.

Now, we describe explicit formula for the action of $R(n_1,n_2)$-matrix. Note that the  $R$-matrix commutes with the action of  $\frak{gl}(N)$ in any representation. Therefore, it is given by a sum of  $\frak{gl}(N)$ invariants (Casimir operators)  with some coefficients. We found that the following description for the ring of the invariants is particularly convenient for the fundamental representations: consider the operators defined by
\be
\label{Omg}
\Omega=\sum\limits_{k=1}^{N}\, \psi_{k}^{\ast} \otimes \psi_{k}, \ \ \Omega^{\ast}=\sum\limits_{k=1}^{N}\, \psi_{k} \otimes \psi_{k}^{\ast}
\ee
Then the operators defined by:
\be
\Xi_{m}=\Omega^m\, {\Omega^{\ast}}^m
\ee
generate the ring of invariants. For example, $\Xi_{0}=1$ and
\be
\Xi_{1}=\Omega\, \Omega^{\ast}=\sum\limits_{i,j} \psi_{i}^{\ast} \psi_{j}\otimes \psi_{i} \psi_{j}^{\ast}=\textrm{Id}-\sum\limits_{i,j}  \psi_{j}\psi_{i}^{\ast}\otimes \psi_{i} \psi_{j}^{\ast}
\ee
and using (\ref{Ephi}) we conclude that;
\be
\Xi_{1}=\textrm{Id}-\textrm{P}, \ \ \textrm{P}=\sum\limits_{i,j=1}^{N} E_{j i} \otimes E_{i j}
\ee
where $\textrm{Id}$ and $\textrm{P}$ are the identity and permutation operator respectively.  The $R$-matrix for the first fundamental representation (\ref{Rsq}) can be expressed through these invariants as follows:
\be
R(u)=\Xi_{0}-\dfrac{\hbar}{u+\hbar} \Xi_{1}
\ee
Computing explicitly the operator $(\ref{Rnm})$ for different choices of $n_1, n_2$ and $N$, and expanding the answer in sum of $\Xi_k$ we found that the
$\frak{gl}(N)$ $R$-matrix for $\bigwedge^{n_{1}} V \otimes \bigwedge^{n_{2}} V$ has the following stable form ( does not depend on $n_1, n_2$ and $N$ ):
\be
\label{Rwed}
\begin{array}{|c|}
\hline \\
 \ \ R(u)=\sum\limits_{k=0}^{\infty} \left(\, \dfrac{(-\hbar)^k}{k! \prod\limits_{m=1}^{k}(u+m \hbar)}\, \right) \,\Xi_{k} \ \ \\
\\
\hline
\end{array}
\ee
Obviously, for $m>\min(n_1,n_2)$ the operators $\Xi_{m}$ vanish on $\bigwedge^{n_{1}} V \otimes \bigwedge^{n_{2}} V$, making the last infinite sum well defined. The stability of this operator allows extend it to the case of $\frak{gl}_{\infty}$.
\subsection{$R$-matrix for the fundamental representations of~$\frak{gl}_{\infty}$ \label{ffs}}
In section~\ref{instR} we have shown that the $R$-matrix for the instanton moduli space can be constructed as certain infinite product of $R$-matrices acting in
\be
\label{ar}
H^{\bullet}_{T} \Big({\cal{A}}_{\infty}(\delta_0+\delta_m) \Big) =H^{\bullet}_{T} \Big({\cal{A}}_{\infty}(\delta_0)  \Big)\otimes H^{\bullet}_{T} \Big({\cal{A}}_{\infty}(\delta_m)  \Big)
\ee
%$$
%{\cal{A}}_{\infty}(\delta_0)=\coprod\limits_{\{n_i\}} {\cal{A}}_{\infty}(n_i,\delta_0)
%$$
In the finite case of $A_{n}$ the space
$H^{\bullet}_T \Big( {\cal{A}}_{n}(\delta_k) \Big)$ has a natural structure of irreducible $\frak{gl}_{n+1}$ module:
$$
H^{\bullet}_T \Big( {\cal{A}}_{n}(\delta_k)\Big)={\bigwedge}^{k} V
$$
where $V$ is the first fundamental representation.

In $\frak{gl}_{\infty}$-case, however, there exist a well known automorphism  acting on the fundamental weights by $\gamma (\delta_{k})=\delta_{k+1}$,
such that all fundamental representations are equivalent to each other, with isomorphism explicitly given by:
$$
\gamma E_{i,j} \gamma^{-1} =  E_{i+1,j+1}
$$
Therefore, $\frak{gl}_{\infty}$-modules $H^{\bullet}_T \Big( {\cal{A}}_{\infty}(\delta_k) \Big)$ are equivalent for all $k$.

In this section we consider the so called half-infinite wedge product $\bigwedge^{\frac{\infty}{2}} V$ for the fundamental representation $V$ of $\frak{gl}_{\infty}$. As $\frak{gl}_{\infty}$-module  $\bigwedge^{\frac{\infty}{2}} V$ splits into the sum of the irreducible representations ${\cal{F}}_{k}$ called \textsl{charge $k$ fermion Fock spaces}:
\be
{\bigwedge}^{\frac{\infty}{2}} V = \bigoplus\limits_{k\, \in\, {\mathbb{Z}}} {\cal{F}}_{k}
\ee
The main reason for considering these modules is the following identification:
\be
\label{indf}
\begin{array}{|c|}
\hline\\
 \ \ \ H^{\bullet}_T \Big({\cal{A}}_{\infty}(\delta_k) \Big)\simeq {\cal{F}}_{k} \ \ \
\\
\\
\hline
\end{array}
\ee
Here we briefly remind this simple construction and describe isomorphisms $\gamma : {\cal{F}}_{k}\rightarrow {\cal{F}}_{k+1}$ of $\frak{gl}_{\infty}$-modules mentioned above.  Details can be found in~\cite{MJD}.

Let $V$ be a space spanned by $e_{i}, \ \ i\in {\mathbb{Z}}+1/2$. Let $\psi_i$ and $\psi_i^{\ast}$ $ i\in {\mathbb{Z}}+1/2$ be the generators of infinite Clifford algebra ${\cal{G}}$ defined by relations (\ref{cl}). Consider the following formal half-infinite wedge product called \textit{vacuum}:
\be
\label{vac}
\textsf{vac}:=e_{1/2}\wedge e_{3/2}\wedge....
\ee
Note that the half of the generators with action defined by (\ref{ac}) annihilate the vacuum:
$$
\psi^{\ast}_{-i}\textsf{vac}= \psi_{i} \textsf{vac} =0 , \ \ \textrm{for} \ \ \ i>0
$$
So that the operators $\psi^{\ast}_{-i}$ and $\psi_{i}$  for $i>0$ are called annihilation operators. The complementary operators
$\psi^{\ast}_{-i}$ and $\psi_{i}$ with $i<0$ are referred to as creation operators. Applied to the vacuum state (\ref{vac}) they "create" a new half-infinite wedge product called  \textit{pure fermion state} $\psi_{n_1}...\psi_{n_r}\psi^{\ast}_{m_1}...\psi^{\ast}_{m_r} \textsf{vac}$. The half-infinite wedge product $\bigwedge^{\frac{\infty}{2}} V$ is defined as the space spanned by the pure fermion states:
\be
{\bigwedge}^{\frac{\infty}{2}} V=\textrm{Span}\{\psi_{n_1}...\psi_{n_r}\psi^{\ast}_{m_1}...\psi^{\ast}_{m_r} \textsf{vac} : n_i, m_j \in {\mathbb{Z}}+1/2 \}
\ee
such that we can also write
\be
{\bigwedge}^{\frac{\infty}{2}} V={\cal{G}} \cdot \textsf{vac}
\ee
The algebra ${\cal{G}}$  and its module $\bigwedge^{\frac{\infty}{2}} V$ carry two important gradings called \textit{charge} and \textit{energy}. These gradings are defined for the generators as follows:
\be
\begin{array}{c|c|c}
 & charge & energy\\
 \hline
 \psi_n &1  & n\\
 \hline
 \psi^{\ast}_n & -1 & n
\end{array}
\ee
using this, the grading of $\bigwedge^{\frac{\infty}{2}} V$ is defined by the charge and energy of the pure states:
\be
\left\{\begin{array}{l}
\textrm{charge  (or energy ) of} \ \ \textsf{vac} =0\\
\\
\textrm{charge  (or energy ) of} \ \  a \cdot\textsf{vac} = \textrm{charge  (or energy ) of} \ \ a
\end{array}\right.
\ee
Denote by ${\cal{F}}^{(d)}_{l}$ the subspace of ${\cal{F}}$ spanned by the pure states with charge $l$ and energy $d$:
\be
{\cal{F}}^{(d)}_{l}=\textrm{Span}\{ \psi_{m_{1}}...\psi_{m_{r}}\psi^{\ast}_{n_{1}}...\psi^{\ast}_{n_{s}} \textsf{vac}\, :\, r-s=l, \, \sum_{i=1}^{r} m_{i} + \sum_{j=1}^{s} n_{j}=d \} \ \
\ee
We  also denote ${\cal{F}}_l=\bigoplus_{d}  {\cal{F}}^{(d)}_{l} $, such that $\bigwedge^{\frac{\infty}{2}} V=\bigoplus_{l}  {\cal{F}}_{l} $.
The generators of $\frak{gl}_{\infty}$ act on this space by normally ordered operators $E_{i,j}=\textrm{:}\psi_{i}\psi_{j}^{\ast}\textrm{:}$. Obviously, these generators have charge zero, therefore each subspace ${\cal{F}}_l$ is $\frak{gl}_{\infty}$-invariant. Moreover, any element of ${\cal{F}}_0$ can be obtained from $\textsf{vac}$ by applying elements from $U(\frak{gl}_{\infty})$, such that:
\be
\label{Uo}
{\cal{F}}_0=U(\frak{gl}_{\infty})\,\textsf{vac}
\ee
Thus, ${\cal{F}}_0$ is an irreducible highest weight module of $\frak{gl}_{\infty}$.
Consider an operator $\gamma \in \textrm{End}(\bigwedge^{\frac{\infty}{2}} V)$  defined by:
\be
\label{gam}
\gamma(e_{i_1}\wedge e_{i_2}\wedge...)=e_{i_1+1}\wedge e_{i_2+1}\wedge...
\ee
Obviously it gives a chain of isomorphisms $\gamma :\, {\cal{F}}_l\rightarrow {\cal{F}}_{l+1}$ for $\frak{gl}_{\infty}$-modules. Such that
all $ {\cal{F}}_l$ are equivalent to $ {\cal{F}}_0$. We can define a \textit{shifted vacuum} by:
\be
\textsf{vac}_{k}=\gamma^{k}\, \textsf{vac} \in {\cal{F}}_k
\ee
such that, analogously to (\ref{Uo}) we have:
\be
{\cal{F}}_k=U(\frak{gl}_{\infty})\,\textsf{vac}_{k}
\ee
Note, that the shifted vacuum $\textsf{vac}_{k}$ has the minimal energy among the pure states ${\cal{F}}_k$ of charge $k$  equal to $k^2/2$.
\subsection{Formula for the instanton $R$-matrix}
Now we are in the position to summarize the results of the whole section and give the formula for the instanton $R$-matrix.
First of all, we have the product formula (\ref{prodFor}). As follows from (\ref{ar}) and identification (\ref{indf}), each term of this infinite product
acts in the tensor product of two Fock spaces with charges shifted by $n$: $R_{n}(u)\in \textrm{End}( {\cal{F}}_{0} \otimes {\cal{F}}_{n})$. The matrix $R_{0}(u)$ can be identified with the matrix $R(u)$  calculated by fusion procedure (\ref{Rwed}). Note that the expression (\ref{Rwed}) was computed for the finite values of $N$, but in the end does not depend on $N$ explicitly and therefore can be extended to the case $N=\infty$.

Using the automorphism $\gamma: {\cal{F}}_{l}\rightarrow {\cal{F}}_{l+1}$ described in the previous section, we can easily construct $R_{n}(u)$ from the zeroth one:
$$
R_{n}(u)=(1\otimes \gamma)^n R_0(u) (1\otimes \gamma)^{-n}
$$
and for the instanton $R$-matrix we obtain:
\be
\label{rr}
{\cal{R}}(u)=\prod\limits_{n\in{\mathbb{Z}}}^{\leftarrow} (1\otimes \gamma)^n\, R(u+n t_{1})\, (1\otimes \gamma)^{-n}
\ee
where the arrow stands for the product ordered from the right to the left. Analogously to (\ref{Omg}) introduce:
\be
\nonumber
\Omega_{n}=(1\otimes \gamma)^n\Big(\sum\limits_{k \in{\mathbb{Z}} +\frac{1}{2}}\, \psi_{k}^{\ast} \otimes \psi_{k}\Big)(1\otimes \gamma)^{-n}=\sum\limits_{k \in{\mathbb{Z}} +\frac{1}{2}}\, \psi_{k}^{\ast} \otimes \psi_{k+n}, \\
\label{Omn} \nonumber
\Omega_{n}^{\ast}=(1\otimes \gamma)^n\Big(\sum\limits_{k \in{\mathbb{Z}} +\frac{1}{2}}\, \psi_{k} \otimes \psi_{k}^{\ast}\Big)(1\otimes \gamma)^{-n}=\sum\limits_{k \in{\mathbb{Z}} +\frac{1}{2}}\, \psi_{k} \otimes \psi_{k+n}^{\ast},
\ee
such that
$$
R_{n}(u)=(1\otimes \gamma)^n\, R(u)\, (1\otimes \gamma)^{-n}=\sum\limits_{k=0}^{\infty} \left(\, \dfrac{(-\hbar)^k}{k! \prod\limits_{m=1}^{k}(u+m \hbar)}\, \right) \,\Omega_{n}^{k} {\Omega_{n}^{\ast}}^{k}
$$
Substituting this and $h=t_1+t_2$ into (\ref{rr}) we obtain the explicit formula for the $R$-matrix:
\be
\label{RF}
\begin{array}{|c|}
\hline\\
\ \ {\cal{R}}(u)=\prod\limits_{n\in{\mathbb{Z}}}^{\leftarrow}\left( \sum\limits_{k=0}^{\infty} \left(\,  \dfrac{(-1)^k (t_1+t_2)^k}{k! \prod\limits_{m=1}^{k}(u+(m+n) t_1+m t_2)}\, \right) \,\Omega_{n}^{k} {\Omega_{n}^{\ast}}^{k}  \right) \ \
\\
\\
\hline
\end{array} \ \
\ee

\section{Expansion of the instanton $R$-matrix \label{expsec}}
In this section we study the Taylor expansion of the instanton $R$-matrix ${\cal{R}}(u)$ and its vacuum matrix element ${\cal{T}}(u)$ in the parameter $u$. We show, that the coefficients of both series can be expressed through bosons (\ref{bosons}) using the boson-fermion correspondence. Moreover, we prove that the expansion of ${\cal{T}}(u)$ completely defines the expansion of ${\cal{R}}(u)$, i.e. the expression in terms of bosons for the $R$-matrix ${\cal{R}}(u)$  can be reconstructed from  ${\cal{T}}(u)$, by doubling the bosons $\alpha_{n}\rightarrow \alpha_{n}\otimes 1 - 1\otimes\alpha_n$.

\subsection{Vacuum matrix element of the instanton $R$-matrix \label{rexpans}}
Let us briefly summarize the results of the previous sections for ${\cal{R}}(u)$. The instanton  $R$-matrix is defined as:
$$
{\cal{R}}(u) = \textrm{Stab}_{-\fC}^{-1}\, \textrm{Stab}_{\fC} \in \textrm{End }\Big( H_{T}^{\bullet} ( \coprod\limits_{n=0}^{\infty} {\cal{M}}^A(n,2)) \Big) \otimes \mathbb{C}(\frak{t})
$$
and the fixed set is a disjoint union of the following components:
$$
{\cal{M}}^A(n,2)=\coprod\limits_{n_1+n_2=n}\, \textrm{Hilb}_{n_1}\times \textrm{Hilb}_{n_2}
$$
In this case $A={\mathbb{C}}^{\ast}$ and is acts on ${\cal{M}}(n,2)$ by elements of the form (\ref{chm}), so there are only two chambers
$\fC=\{u>0\}$ and $-\fC$. The ordering $\succeq$ defined on the components of ${\cal{M}}^A(n,2)$ by the chamber $\fC$ is of the form:
$$
\textrm{Hilb}_{n_1}\times \textrm{Hilb}_{n_2} \succeq \textrm{Hilb}_{m_1}\times \textrm{Hilb}_{m_2}\ \ \Leftrightarrow \ \ n_1\leq m_1
$$
We denote $Z_{i}=\textrm{Hilb}_{n-i}\times \textrm{Hilb}_{i}$ and such that the ordering on the fixed set is:
$$
Z_0\preceq Z_1\preceq...\preceq Z_n
$$
Note, that $\textrm{Hilb}_{0}\simeq \textrm{pt}$, thus  $\textrm{Hilb}_{n}\times\textrm{pt}\simeq \textrm{Hilb}_{n}$ and  $\textrm{pt}\times  \textrm{Hilb}_{n} \simeq \textrm{Hilb}_{n}$ are the minimal and maximal components of the fixed set.  Now, let
$$
i: {\cal{M}}^A(n,2)\hookrightarrow {\cal{M}}(n,2)
$$
be the inclusion of the fixed set, then the corresponding map:
$$
i^{\ast}: H_{T}^{\bullet}\Big({\cal{M}}(n,2)\Big)\otimes \mathbb{C}(\frak{t}) \longrightarrow H_{T}^{\bullet}\Big({\cal{M}}^A(n,2)\Big)\otimes \mathbb{C}(\frak{t})
$$
becomes an isomorphism after localization and we can rewrite $R$-matrix as:
$$
{\cal{R}}(u) = \Big(i^{\ast}\,\textrm{Stab}_{-\fC}\Big)^{-1}\, \Big(i^{\ast}\, \textrm{Stab}_{\fC} \Big)
$$
By the first property of the stable map (section \ref{se}) we have $\textrm{Supp} \Big(\textrm{Stab}_{\fC}(Z) \Big) \subset \, \textrm{Slope}(Z) $ for any component
$Z \in {\cal{M}}^A(n,2)$. It means that
$$
i^{\ast}_{Z_j} \textrm{Stab}_{\fC}(Z_i)=0 \ \  \textrm{for} \ \  Z_j \succeq Z_i
$$
where by $i^{\ast}_{Z_j}$ we denote the restriction of a cohomology class to the component $Z_j$. Therefore, we conclude that in the ordered  basis $Z_i$, the corresponding maps are given by some upper and lower block-triangular matrices:
$$
i^{\ast} \textrm{Stab}_{\fC}=\left(\begin{array}{cccc}
S_{0 0}&S_{0 1}&S_{0 2}&..\\
0&S_{1 1}&S_{1 2}&..\\
0&0&S_{2 2}&..\\
..&..&..&..
\end{array}\right)\,
i^{\ast} \textrm{Stab}_{-\fC}=\left(\begin{array}{cccc}
U_{0 0}&0&0&..\\
U_{1 0}&U_{1 1}&0&..\\
U_{2 0}&U_{2 1}&U_{2 2}&..\\
..&..&..&..
\end{array}\right)
$$
The block submatrices $S_{i j}$ give maps:
$$
S_{i j}: H^{\bullet}_{T}{( Z_j )}\longrightarrow H^{\bullet}_{T}{( Z_i )}
$$
and similarly for $U_{i j}$. From the triangular form of the matrices we obtain:
\be
\label{vacelem}
{\cal{T}}(u)={\cal{R}}(u)_{0 0}=U_{0 0}^{-1} S_{0 0} \in \textrm{End} \left(H^{\bullet}_{T}{\Big( \coprod\limits_{n=0}^{\infty} \textrm{Hilb}_{n} \Big)}\right) \otimes \mathbb{C}(\frak{t})
\ee
By the second property of the stable map we have:
\be
S_{i i}=e(N[i]_{-}), \ \ U_{i i}=e(N[i]_{+})
\ee
where $N[i]$ is the normal bundle to the component $Z_{i}$ in ${\cal{M}}(n,2)$. Therefore, using (\ref{polar})  we obtain:
$$
S_{0 0}=e(N_-), \ \ U_{0 0}=e(N_{+})=e(N_{-}\otimes \hbar)
$$
where $N_-$ is the component of the normal bundle to the minimal component $\textrm{Hilb}_n$ in ${\cal{M}}(n,2)$. From section \ref{tb} we have:
\be
N_-={\cal{V}}\otimes u
\ee
Therefore the operator ${\cal{T}}(u)$ acts in ${\cal{F}}$ as the following characteristic class:
\be
\label{mec}
\begin{array}{|c|}
\hline\\
{\cal{T}}(u)=\dfrac{e\Big({\cal{V}} \otimes u\Big)}{e\Big({\cal{V}} \otimes u \otimes \hbar\Big)}\\
\\
\hline
\end{array}
\ee
Remind that the fixed set $\textrm{Hilb}^{T}_{n}$ is the disjoint union of points, whose number is equal to the number of partitions of $n$. Each point is naturally labeled by a partition $\lambda$ with $|\lambda|=n$.  The classes of the fixed points $[\lambda]\in H_{T}^{\bullet}(\textrm{Hilb}_{n})$ correspond to the Jack polynomials in the Fock space $j_{\lambda} \in {\cal{F}}$. From (\ref{Nekf}), the restriction of the equivariant Euler class to the point $\lambda$ has the form:
\be
\left.e(N_{-}\otimes u)\right|_{[\lambda]}=\prod\limits_{(i,j)\in \lambda}\,\Big(u+(i-1)t_1 +(j-1) t_2 \Big)
\ee
Thus, taking into account that $\hbar=t_1+t_2$ we obtain that the vacuum matrix element ${\cal{T}}(u)$  acts in the basis of Jack polynomials diagonally with the following eigenvalues:
\be
\label{ev}
{\cal{T}}(u) \, j_{\lambda}=\dfrac{\prod\limits_{(i,j)\in \lambda}\,\Big(u+(i-1)t_1 +(j-1) t_2 \Big)}{\prod\limits_{(i,j)\in \lambda}\,\Big(u+i t_1+j t_2 \Big)} j_{\lambda}
\ee
As explained in the introduction, (\ref{mec}) can be considered as the generating function for multiplication by Chern classes of tautological bundle
$c_{n}({\cal{V}})$:
\be
\label{expan}
{\cal{T}}(u)=1-\dfrac{\hbar\, \textrm{rk}  {\cal{V}} }{u}+\dfrac{\hbar c_1 ({\cal{V}}) + \hbar^2 \textrm{rk}  {\cal{V}}(\textrm{rk}  {\cal{V}}+1)/2}{u^2} +O(u^{-2})
\ee
Expanding (\ref{ev}) into the Taylor series in $u^{-k}$, and matching terms with one in (\ref{expan}) we obtain:
\be
\label{egv}
\textrm{rk} \,j_{\lambda} =|\lambda| \, j_{\lambda}, \ \ \ c_1 ({\cal{V}}) j_{\lambda}= \Big(\sum\limits_{(i,j)\in \lambda}\, (i-1)t_1 +(j-1) t_2 \,\Big)\,j_{\lambda}
\ee
The matrix elements ${\cal{T}}(u)$, in principle, can be computed explicitly from (\ref{RF}), what gives explicit expression for the operators $c_{n}({\cal{V}})$ in terms of fermions $\psi_n$. To pass to the Lehn-like formulae (\ref{lf}) for the Chern classes, we  rewrite the fermions through the bosons $\alpha_n$ which is, of course, the well known boson-fermion correspondence.

\subsection{The boson-fermion correspondence}
Let us consider the charge zero fermion Fock space ${\cal{F}}_{0}$. As we discuss in section \ref{ffs} this space is spanned by the pure fermion states of charge zero and is graded by the energy:
$$
{\cal{F}}_{0}=\bigoplus\limits_{d=0}^{\infty} {\cal{F}}^{(d)}_{0}
$$
Recall, that there is a natural labeling of the pure states by the Maya diagrams or, equivalently, by the partitions \cite{MJD}.  Each fermion pure state $v_\lambda\in {\cal{F}}^{(d)}_{0}$ corresponds to certain partition $\lambda$ with $|\lambda|=d$.

Let ${\cal{F}}\simeq\mathbb{C}[p_1,p_2,...]$ be the boson Fock space described in the introduction.  If we define $\deg p_n =n$, then:
$$
{\cal{F}}=\bigoplus\limits_{d=0}^{\infty} {\cal{F}}^{(d)}
$$
where ${\cal{F}}^{(d)}$ is the space of polynomials of degree $d$. The dimension of ${\cal{F}}^{(d)}$ is equal to the number of partitions of $d$. The natural basis of ${\cal{F}}^{(d)}$ labeled by partitions $\lambda$ is given by the Schur polynomials, which can be defined explicitly as:
$$
s_{\lambda}=\det_{ij} g_{\lambda_i-i+j}
$$
where $g_{k}$ are the coefficients of the expansion:
$$
\exp\Big( \sum\limits_{n=1}^{\infty} \dfrac{p_n}{n} z^n \Big)=\sum\limits_{k=0}^{\infty}\,g_k z^k
$$
$\deg s_\lambda =|\lambda|$ such that $s_\lambda\in {\cal{F}}^{(|\lambda|)}$.

The boson-fermion correspondence is the isomorphism of the fermion Fock  space ${\cal{F}}_{0}$ and the boson Fock space ${\cal{F}}$ induced by the identification of their basis elements:
\be
\label{iden}
v_\lambda \rightarrow s_\lambda
\ee
This identification leads to certain relations between bosons (\ref{bosons}) and fermions (\ref{ac}), such that they can be expressed through each other. The action of the bosons can be expressed through the fermions as:
\be\label{breal}
\alpha_{n}=\sum\limits_{m\in \mathbb{Z}+\frac{1}{2}} :\psi_i \psi_{i+n}^{\ast}:
\ee
To show this, its enough to check that the operator $a_{-n}$ creates a state $p_n$ in ${\cal{F}}=\mathbb{C}[p_1,p_2,...]$ under the identification (\ref{iden}), and the operator $\alpha_n$ kills the vacuum.

Note that the operators (\ref{breal}) are neutral i.e. they do not change the charge:
$$
a_{n} : {\cal{F}}_l \rightarrow {\cal{F}}_l
$$
 Thus, each charge subspace of the fermion Fock space is identified with a copy of the boson space ${\cal{F}}_l\simeq {\cal{F}}$.  The formulae for the correspondence that treat all these spaces at one go, using generating functions can be found in \cite{MJD}.  Using (\ref{breal}), we can rewrite certain neutral infinite fermion sums in terms of bosons. In particular, we can rewrite the matrix elements of the instanton $R$-matrix  (\ref{RF}) directly through the bosons.
\subsection{Expansion of the vacuum matrix element ${\cal{T}}(u)$}
The instanton $R$-matrix (\ref{RF}) is neutral, thus defines an operator:
\be
{\cal{R}}(u): {\cal{F}}_0 \otimes {\cal{F}}_0 \rightarrow  {\cal{F}}_0 \otimes {\cal{F}}_0
\ee
The operator ${\cal{T}}(u)$ is the matrix element of the $R$-matrix for the subspace $\textsf{vac} \otimes {\cal{F}}_0 \subset  {\cal{F}}_0 \otimes {\cal{F}}_0$. The Taylor expansion of the vacuum element has the form
\be
\label{evme}
{\cal{T}}(u)=1+\sum\limits_{n=1}^{\infty} \, u^{-n} \, {\cal{H}}_{n}
\ee
From explicit formula for $R$-matrix (\ref{RF}) and correspondence (\ref{breal}) we obtain expressions for the first two coefficients:
\begin{small}
$$
\label{Tc}
\begin{array}{l}
{\cal{H}}_{1}=\hbar \sum\limits_{n=1}^{\infty} \alpha_{-n}\alpha_{n},\\
\\
{\cal{H}}_{2}=\dfrac{\hbar}{2} \sum\limits_{m,n=1}^{\infty}\Big(  t_1 t_2 \alpha_{-m} \alpha_{-n} \alpha_{n+m}-\alpha_{-m-n} \alpha_{m} \alpha_n\Big)
+\dfrac{\hbar^2}{2} \Big(\sum\limits_{n=1}^{\infty} \alpha_{-n}\alpha_n\Big)^{2} +\dfrac{\hbar^2}{2} \sum\limits_{n=1}^{\infty} n \alpha_{-n}\alpha_n
\end{array}
$$
\end{small}
In the same way we can obtain explicit formulae for higher coefficients, but their form is more complicated.  Now, comparing this expansion with (\ref{expan}) we obtain:
\be
\label{rank}
\textrm{rk} ({\cal{V}})=\sum\limits_{n=1}^{\infty} \alpha_{-n}\alpha_{n}=\sum\limits_{n=1}^{\infty} n p_n \dfrac{\partial}{\partial p_{n}}
\ee
Note that the rank of the tautological bundle ${\cal{V}}$ over $\textrm{Hilb}_m$ is equal $\textrm{rk} ({\cal{V}})~=~m$. In full agreement, the operator (\ref{rank}), obviously, acts by multiplication on $m$ in $H^{\bullet}_{T}(\textrm{Hilb}_m)$.

Comparing $u^{-2}$ term with (\ref{expan}) we obtain an explicit expression for the first Chern class:
\begin{small}
\be
\begin{array}{r}
c_{1}({\cal{V}})=\dfrac{1}{2}\sum\limits_{m,n=1}^{\infty}\Big( t_1 t_2 \alpha_{-m} \alpha_{-n} \alpha_{n+m} - \alpha_{-m-n} \alpha_{n} \alpha_{m}\Big)+
\dfrac{\hbar}{2}\, \sum\limits_{n=1}^{\infty}\,( n-1) \alpha_{-n}\alpha_{n}
\end{array} \ \
\ee
\end{small}
which is in full agreement with the result of Lehn.   The eigenvalues of this operator in the basis of
the Jack polynomials are given by (\ref{egv}).
\subsection{Expansion of  ${\cal{R}}(u)$}
Let us consider the tensor square of the boson Fock space, and the operators acting  as follows:
\be
\alpha^{\pm}_{n}=\alpha_{n} \otimes 1\pm 1\otimes  \alpha_{n}
\ee
The commutation relations for these operators are:
\be
[\alpha^{+}_{n},\alpha^{+}_{m}]=2 \delta_{n+m}, \ \ [\alpha^{-}_{n},\alpha^{-}_{m}]=2 \delta_{n+m}, \ \ [\alpha^{+}_{n},\alpha^{-}_{m}]=0
\ee
The operators $\alpha^{\pm}_{n}$ generate a couple of Heisenberg algebras $\frak{h}^{\pm}$ such that the elements of
$\frak{h}^{+}$ commute with elements of $\frak{h}^{-}$. We also denote by ${\cal{F}}^{\pm}$ the Fock representations for $\frak{h}^{\pm}$ respectively.
Thus we can write:
$$
{\cal{F}} \otimes {\cal{F}} = {\cal{F}}^{+} \otimes {\cal{F}}^{-}
$$

\noindent
\textbf{Proposition:} The instanton $R$-matrix ${\cal{R}}(u)$ is an element of $\frak{h}^{-}$.

\vspace{1mm}
\noindent
\textit{Proof}  By construction, ${\cal{R}}(u)$ is given by a product (\ref{prodFor}) of rational $R$-matrices corresponding to the half-infinite wedge product representation $\Lambda^{\frac{\infty}{2}} V$ for $\frak{gl}_{\infty}$. The general property of the rational $R$-matrix is that it commutes with the action of $\frak{gl}_{\infty}$.Thus, each term in (\ref{prodFor}), and therefore ${\cal{R}}(u)$ itself commute with the action of $\frak{gl}_{\infty}$. It implies, in particular that:
$$
[{\cal{R}}(u), a^{+}_{n}]=0, \ \ \forall \, n
$$
now, the proposition follows from the fact that $\frak{h}^{\pm}$ acts on ${\cal{F}}^{\pm}$ irreducibly $\Box$.

Therefore, the coefficients of expansion of the instanton $R$-matrix can be expressed through $\alpha^{-}_{n}$. In fact, to obtain the expansion of $R$-matrix it is enough to know the expansion of its vacuum matrix element. Indeed, denote by $A_{\emptyset}$ the vacuum matrix element of an operator $A \in \textrm{End} ({\cal{F}}^{2})$.

\noindent
\textbf{Proposition:} An element $A \in {\frak{h}}^{-}$ is defined uniquely by is matrix element~$A_{\emptyset}$. If $A_{\emptyset}$ is already expressed through the normally ordered combinations of the bosons $\alpha_n$, then $A$ is obtained from $A_{\emptyset}$ by the change of variables $\alpha_n\rightarrow \alpha_n^{-}$.

\smallskip

\noindent
\textit{Proof.}
% Let  $I=\{n_{i}\}$ be a set of nonnegative numbers $n_{i}\geq 0$ labeled by  ${i\in \mathbb{Z}}$ such that only finitely many $n_{i}\neq 0$.
%Consider normally ordered monomial $a^{-}_{I}\in End^{-}$ defined by:
%\be
%a^{-}_{I}=:\prod\limits_{i\in \mathbb{Z}}  (a^{-}_{i})^{n_{i}}:
%\ee
Let $\lambda=\lambda_1\geq ...\geq \lambda_m$ be a partition. Consider the normally ordered elements $:\alpha_{\lambda}^{-}:=\alpha_{\lambda_1}^{-}...\alpha_{\lambda_m}^{-}$ forming a basis in $\frak{h}^{-}$.
%Elements $a^{-}_{I}$ form a basis of $End^{-}$.
We have
$$
\alpha^{-}_{\lambda}=\alpha_{\lambda} \otimes 1 +...+  1 \otimes \alpha_{\lambda}
$$
where dots stand for the mixed terms of the form $\alpha_{\mu}\otimes \alpha_{\nu}$.
It is clear that the vacuum matrix element is defined by the first term $(:a^{-}_{\lambda}:)_{\emptyset}=a_{\lambda}$. Note, that $a_{\lambda}$ is also normally ordered. Moreover, given a normally ordered monomial $a_{\lambda}$ then $A=a^{-}_{\lambda}$ is a unique element in ${\frak{h}}^{-}$ with the property $A_{\emptyset}=\alpha_{\lambda}$. $\square$

For example, we can find the first two coefficients of the expansion
$$
{\cal{R}}(u)=1+\sum\limits_{n=1}^{\infty}\, u^{-n}\,  {\cal{R}}_{n}
$$
from the known coefficients of ${\cal{T}}(u)$ (\ref{evme}). We obtain:
$$
\begin{array}{l}
{\cal{R}}_{1}=\hbar \sum\limits_{n=1}^{\infty}\, \alpha_{-n}^{-} \alpha_{n}^{-},\\
\\
{\cal{R}}_{2}=\dfrac{\hbar}{2} \sum\limits_{m,n=1}^{\infty}\Big(\, t_1 t_2 \, \alpha_{-n}^{-} \alpha_{-m}^{-} \alpha_{n+m}^{-} - \, \alpha_{-n-m}^{-} \alpha_{m}^{-} \alpha_{n}^{-}\Big) + \dfrac{\hbar^2}{2} \sum\limits_{m,n=1}^{\infty}\, \alpha_{-n}^{-}\alpha_{-m}^{-}\alpha_{n}^{-}\alpha_{m}^{-}
\end{array}
$$
To prove it, by the last proposition it is enough to check that:
$$
\Big({\cal{R}}_{1}\Big)_{\emptyset}={\cal{H}}_{1}, \ \ \Big({\cal{R}}_{2}\Big)_{\emptyset}={\cal{H}}_{2}
$$
The last two expressions for ${\cal{R}}_{1}$ and ${\cal{R}}_{2}$ suggest that
there is a wise choice for normalization of the boson operators. Indeed consider the scaled bosons, corresponding to Nakajima operators for equivariant lines in ${\mathbb{C}}^2$:
\be
\textbf{a}_{- n}=- \frac{1}{t_1 t_2} \, \alpha_{-n}^{-}, \ \ \textbf{a}_{n}=  \alpha_{n}^{-}
\ee
and consider the corresponding bosonic field:
\be
{\cal{A}}(z)=\sum\limits_{n =-\infty}^{\infty} \textbf{a}_{n} \,z^{n}
\ee
If we introduce the following notations for the residue of the formal series:
$$
\oint f(z) dz=f_{0}, \ \ \ \textrm{for} \ \ \ f(z)=\sum\limits_{n =-\infty}^{\infty} f_{n} \,z^{n}
$$
then we can rewrite the above formulae for the coefficients of intanton $R$-matrix in the following compact form
\be
\begin{array}{l}
{\cal{R}}_{1}=-\dfrac{\hbar}{t_1 t_2}\oint {\cal{A}}^2(z) dz\\
\\
{\cal{R}}_{2}=\dfrac{\hbar}{2 t_1^2 t_2^2} \oint {\cal{A}}^3(z) dz  +\dfrac{1}{2} \Big(\dfrac{\hbar}{t_1 t_2} \oint {\cal{A}}^2(z) dz\Big)^2
\end{array}
\ee
In particular the expansion of $R$-matrix takes a simple form:
\be
\log {\cal{R}}(u)=1 - u^{-1}\,\dfrac{\hbar}{t_1 t_2}\oint {\cal{A}}^2(z) dz + u^{-2}\, \dfrac{\hbar^2}{2 t_1^2 t_2^2}\oint {\cal{A}}^3(z) dz + ...
\ee
We believe that there is a simple description of the higher coefficients for the $R$-matrix.We leave the analysis  of the coefficients of the instanton $R$-matrix and, hopefully, the explicit combinatorial formulae for them for the future work.


\begin{thebibliography}{12}
\bibitem{OM} A. Okounkov, D. Maulik, Quantum Groups and Quantum Cohomology, to appear.
\bibitem{Ok1}  A. Okounkov, R. Pandharipande, Quantum cohomology of the Hilbert scheme of points in the plane, 	Invent. math. 179: 523–557 arXiv:math/0411210v2
\bibitem{Lehn} M. Lehn, Lectures on Hilbert Schemes, CRM proceedings and lecture notes, Vol. 38, p. 1-30.
\bibitem{Lehn2} M. Lehn,  Chern classes of tautological sheaves on Hilbert schemes of
points on surfaces, Invent. Math. 136 (1999), no. 1, 157-207.
\bibitem{HL} D. Huybrechts, M. Lehn, The Geometry of Moduli Spaces of Sheaves, Cambridge University Press, Sec. Ed. (2010).
\bibitem{Nak1} H.Nakajima, Lectures on Hilbert Schemes of points on surfaces, AMS~(1999)
\bibitem{NakAle} H.Nakajima, Instantons on ALE spaces, quiver varieties, and Kac-Moody algebras, Duke Math. J. Volume 76, Number 2 (1994), 365-416.
\bibitem{NakQuiv} H.Nakajima, Quiver varieties and Kac-Moody algebras, Duke Math. J. Volume 91, Number 3 (1998), 515-560.
\bibitem{Nak4} H. Nakajima, Jack polynomials and Hilbert schemes of points on surfaces, alg-geom/9610021.
\bibitem{Nak5} H. Nakajima, Heisenberg algebra and Hilbert schemes of points on projective surfaces, Ann. of Math. (2) 145 (1997), no. 2, 379–388.

\bibitem{Nak6}  H. Nakajima, Quiver varieties and finite dimensional representations of quantum affine algebras, JAMS, Volume 14, p 145-238
\bibitem{Varagnolo} M. Varagnolo, E. Vasserot, Canonical bases and quiver varieties, Rep. theor. JAMS,
Volume 7, p. 227-258
\bibitem{Michael} M. McBreen, unpublished.

\bibitem{PandLect} Kentaro Hori, Sheldon Katz, Albrecht Klemm, Rahul Pandharipande, Richard Thomas, Cumrun Vafa, Ravi Vakil, Eric Zaslow, Mirror Symmetry,
American Mathematical Society 2003.



\bibitem{ADHM} Michael Atiyah, Vladimir G. Drinfel'd, Nigel. J. Hitchin, Yuri I. Manin, Construction of Instantons, Phys. Lett. A65 (1978) 185-187
\bibitem{CarOk} Erik Carlsson and Andrei Okounkov, Exts and vertex operators, Duke Math. J. Volume 161, Number 9 (2012), 1797-1815, arxiv: 0801.2565
\bibitem{Goem1} V. Baranovsky, Moduli of sheaves on surfaces and action of the oscillator
algebra, J. Differential Geom. 55 (2000), no. 2, 193–227.
\bibitem{Geom2}  K. Costello and I. Grojnowski, Hilbert schemes, Hecke algebras and the
Calogero-Sutherland system, arXiv:math/0310189.
\bibitem{Geom3} I. Grojnowski, Instantons and affine algebras, I. The Hilbert scheme and
vertex operators, Math. Res. Lett. 3 (1996), no. 2, 275–291
\bibitem{Geom4} M. Haiman, Macdonald polynomials and geometry, New Perspectives in
Geometric Combinatorics, MSRI Publications 37 (1999), 207–254.
\bibitem{Geom5} O. Schiffmann and E. Vasserot, The elliptic Hall algebra and the equivariant K-theory of the Hilbert scheme of ${\mathbb{A}}^2$
, arXiv:0905.2555.
\bibitem{MMN1}    A.Mironov, A.Morozov, S.Natanzon, Complete Set of Cut-and-Join Operators in Hurwitz-Kontsevich Theory, Theor.Math.Phys.166:1-22,2011, 	 arXiv:0904.4227
\bibitem{MMN2} A.Mironov, A.Morozov, S.Natanzon, Algebra of differential operators associated with Young diagrams, Journal of Geometry and Physics 62 (2012), pp. 148-155, 	arXiv:1012.0433v1

\bibitem{StrThFirst} E.Witten, "Chern-Simons gauge theory as a string theory",
Prog.Math.133 (1995) 637-678, arXiv: 9207094
\bibitem{StrTh2} M.Aganagic, M.Marino and C.Vafa, "All loop topological string amplitudes from
Chern-Simons theory", Commun.Math.Phys. 247 (2004) 467-512, arXiv: 0206164
\bibitem{StrThLast}
 M. Aganagic, A. Klemm, M. Marino, C. Vafa, "The Topological vertex," Commun.
Math. Phys. 254, 425-478 (2005),  arXiv: 0305132.


\bibitem{SUSY1} N. Nekrasov, Seiberg-Witten prepotential from instanton counting, Adv.
Theor. Math. Phys. 7 (2003), no. 5, 831–864.
\bibitem{SUSY2} N. Nekrasov and S. Shatashvili, Supersymmetric vacua and Bethe
ansatz, Nuclear Phys. B Proc. Suppl. 192/193 (2009), 91–112.
\bibitem{SUSY3} N. Nekrasov and S. Shatashvili, Quantum integrability and supersymmetric vacua, arXiv:0901.4748.
\bibitem{SUSY4} N. Nekrasov and S. Shatashvili, Quantization of integrable systems and
four dimensional gauge theories, XVIth International Congress on Mathematical Physics, 265–289, World Sci. Publ., 2010
\bibitem{SUSYi} A.Mironov, A.Morozov, Sh.Shakirov, Matrix Model Conjecture for Exact BS Periods and Nekrasov Functions, 	JHEP 1002:030,2010, 	 arXiv:0911.5721
\bibitem{SUSY5} L. Alday, D. Gaiotto, and Y. Tachikawa, Liouville correlation functions from four-dimensional gauge theories, Lett. Math. Phys. 91 (2010),
no. 2, 167–197.



\bibitem{Knots1} E. Witten, Quantum Field Theory and the Jones Polynomial, Commun. Math. Phys. 121
(1989) 351
\bibitem{KnotsR} N. Reshetikhin, V. G. Turaev, Invariants of 3-manifolds via link polynomials and quantum groups,  Inventiones mathematicae, Volume 103, Issue 1, pp. 547-597, (1991)
\bibitem{Knots2} M. Aganagic and S. Shakirov, Knot Homology from Refined Chern-Simons Theory,
arXiv:1105.5117
\bibitem{Knots3} P. Dunin-Barkowski, A. Mironov, A. Morozov, A. Sleptsov and A. Smirnov, Superpolynomials for toric knots from evolution induced by cut-and-join operators, arXiv:1106.4305
\bibitem{Knots4} A.Anokhina, A.Mironov, A.Morozov, An.Morozov, Knot polynomials in the first non-symmetric representation, 	arXiv:1211.6375
\bibitem{Knots5} E. Gorsky, q,t-Catalan numbers and knot homology, arXiv:1003.0916
\bibitem{Knots6} D.Galakhov, A.Mironov, A.Morozov, A.Smirnov, On 3d extensions of AGT relation, Theor.Math.Phys. 172 (2012) 939-962; Teor.Mat.Fiz. 172 (2012) 73-99	arXiv:1104.2589
\bibitem{Knots7} L.Faddeev and R.Kashaev, Mod.Phys.Lett. A9 (1994) 427-434, hep-th/9310070;
R. Kashaev, Lett.Math.Phys. 43 (1998) 105-115

\bibitem{QISMfirst} L. Faddeev, How the algebraic Bethe ansatz works for integrable models, Sym´etries quantiques (Les Houches, 1995), 149–219, North-Holland,
1998
\bibitem{QISM1} L. Faddeev, N. Reshetikhin, L. Takhtadzhyan, Quantization of Lie
groups and Lie algebras, Leningrad Math. J. 1 (1990), no. 1, 193–225.
\bibitem{QISM2} M. Jimbo, T. Miwa, Algebraic analysis of solvable lattice models, American Mathematical Society, 1995.
\bibitem{Sklyanin} E. K. Sklyanin,  Quantum Inverse Scattering Method. Selected Topics. 	arXiv:hep-th/9211111v1
\bibitem{MJD} T.Miwa, M. Jimbo, E. Date, Solitons, Cambridge university press, 2000.
\bibitem{Zabrodin} A. Zabrodin, Discrete Hirota's equation in quantum integrable models, 	arXiv:hep-th/9610039
\bibitem{QISMlast} A. Mironov, A. Morozov, Y. Zenkevich, A. Zotov,  Spectral Duality in Integrable Systems from AGT Conjecture, 	arXiv:1204.0913
\bibitem{Brion} M. Brion, Equivariant cohomology and equivariant intersection theory,  Proc. of the NATO Advanced Study Institute on representation theories and algebraic geometry, vol.514 (1997)









\end{thebibliography}
\end{document}